\documentclass[reqno]{amsart}

\usepackage{amssymb,amscd,amsxtra,epic,eepic}

\renewcommand\theenumi{\roman{enumi}}

\numberwithin{equation}{section}
\newcommand\note[1]%
{$^\dagger$\marginpar{\sf\footnotesize{$^\dagger${#1}}}}

\swapnumbers
\newtheorem{theorem}{Theorem}[subsection]
\newtheorem{proposition}[theorem]{Proposition}
\newtheorem{lemma}[theorem]{Lemma}
\newtheorem{corollary}[theorem]{Corollary}

\theoremstyle{definition}
\newtheorem{definition}[theorem]{Definition}
\newtheorem{example}[theorem]{Example}
\newtheorem{remark}[theorem]{Remark}

\newcommand\eu{\mathfrak}
\newcommand\lie{\mathfrak}
\renewcommand\k{\lie{k}} 
\renewcommand\t{\lie{t}}
\newcommand\g{\lie{g}}
\newcommand\h{\lie{h}} 
\newcommand\s{\lie{s}}

\newcommand\bb{\mathbb}

\newcommand\Z{\bb{Z}} 

\newcommand\R{\bb{R}} 
\newcommand\C{\bb{C}}

\newcommand\ca{\mathcal}

\newcommand\hull{\operatorname{hull}}
\newcommand\cone{\operatorname{cone}}

\newcommand\Hom{\operatorname{Hom}}

\newcommand\id{\operatorname{id}}

\newcommand\Lie{\operatorname{Lie}} 

\newcommand\SL{\operatorname{\mathbf{SL}}}
\newcommand\SU{\operatorname{\mathbf{SU}}}
\newcommand\U{\operatorname{\mathbf{U}}}

\newcommand\G{\operatorname{\mathbf{G}}}

\newcommand\pre{\preccurlyeq}
\newcommand\suc{\succcurlyeq}

\newcommand\qu{/\kern-.7ex/}
\newcommand\bigqu{\big/\kern-.85ex\big/}

\newcommand\longhookrightarrow{\lhook\joinrel\longrightarrow}

\newcommand\inj{\longhookrightarrow}
\newcommand\sur{\longrightarrow\kern-1.9ex\to}
\newcommand\iso{\longhookrightarrow\kern-1.9ex\to}

\newcommand\bu{{\scriptscriptstyle\bullet}}
\newcommand\ti{\tilde}
\newcommand\tlambda{{\ti\lambda}}
\newcommand\tmu{{\ti\mu}}
\newcommand\tW{W\kern-0.8em\widetilde{\phantom I}\kern0.3em}
\newcommand\barW{W\kern-0.8em\overline{\phantom I}\kern0.3em}

\newcommand\antiddots{\mathinner{%
\mkern1mu\raise1pt\vbox{\kern7pt\hbox{.}}
\mkern2mu\raise4pt\hbox{.}\mkern2mu\raise7pt\hbox{.}\mkern1mu}}

\newcommand\inv{^{-1}} 

\renewcommand\subset{\subseteq}
\renewcommand\supset{\supseteq}

\newcommand\sst{^{\mathrm{ss}}}
\newcommand\rel{_{\mathrm{rel}}}
\newcommand\com{_{\mathrm{com}}}

\title[Orbits, polytopes, and the Hilbert-Mumford criterion]{Coadjoint
orbits, moment polytopes, and the Hilbert-Mumford criterion}

\author{Arkady Berenstein}

\address{Department of Mathematics, Cornell University, Ithaca, New
York 14853-4201} 

\curraddr{Department of Mathematics, Harvard University, Cambridge,
Massachusetts 02138-2901}

\email{arkadiy@math.harvard.edu}

\author{Reyer Sjamaar}

\address{Department of Mathematics, Cornell University, Ithaca, New
York 14853-4201} 

\thanks{The second author was partially supported by an Alfred
P. Sloan Research Fellowship and by NSF Grant DMS-9703947}

\email{sjamaar@math.cornell.edu}

\date{19 October 1998.  Revised 21 November 1999}

\begin{document} 


\begin{abstract}
Consider a compact Lie group and a closed subgroup.  Generalizing a
result of Klyachko, we give a necessary and sufficient criterion for a
coadjoint orbit of the subgroup to be contained in the projection of a
given coadjoint orbit of the ambient group.  The criterion is couched
in terms of the ``relative'' Schubert calculus of the flag varieties
of the two groups.
\end{abstract}


\maketitle

\tableofcontents

\section{Introduction}

Let $K$ be a compact connected Lie group and let $\ti K$ be a closed
connected subgroup.  Let $f$ denote the inclusion of $\ti K$ into $K$,
$f_*\colon\ti\k\to\k$ the induced embedding of Lie algebras, and
$f^*\colon\k^*\to\ti\k^*$ the dual projection.  In this paper we
present a solution to the following problem.

\begin{enumerate}
\item\label{item;geometric}
Let $\ca O$ be a coadjoint orbit of $K$.  Its projection $f^*(\ca O)$
is a $\ti K$-stable subset of $\ti\k^*$.  Which coadjoint orbits of
$\ti K$ are contained in $f^*(\ca O)$?
\end{enumerate}

Select maximal tori $T$ in $K$ and $\ti T$ in $\ti K$, and Weyl
chambers $\t^*_+$ in $\t^*$ and $\ti\t^*_+$ in $\ti\t^*$, where $\t$
and $\ti\t$ denote the Lie algebras of $T$, resp.\ $\ti T$.  The set
$f^*(\ca O)$ is completely determined by the intersection $\Delta(\ca
O) =f^*(\ca O)\cap\ti\t^*_+$, and according to a result of Heckman
\cite{heckman;projections;inventiones}, $\Delta(\ca O)$ is a convex
polytope.  It can be described as the support of the asymptotic
multiplicity function, for which Heckman found an explicit formula.
However, because the formula involves an alternating sum in which many
cancellations take place, it is in practice impossible to determine
the polytope from it.  Our goal is to write down as explicitly as
possible the inequalities defining the polytope.

Our results are inspired by a remarkable paper of Klyachko
\cite{klyachko;stable-bundles-hermitian}, in which he solved the case
of the diagonal embedding of $\U(n)$ into $\U(n)\times\U(n)$.  See
also Fulton's recent survey paper \cite{fulton;eigenvalues-sums}.  As
in \cite{klyachko;stable-bundles-hermitian}, the proof of our
inequalities relies on the Hilbert-Mumford criterion for semistability
of flags with respect to equivariant ample line bundles.  A related
result describing the ``general'' faces of the polytope was arrived at
independently by Brion \cite{brion;general-faces}.

Recall that the restriction of $f^*$ to $\ca O$ is nothing but the
moment map for the action of $\ti K$ on the homogeneous symplectic $
K$-manifold $\ca O$ and that $\Delta(\ca O)$ is its moment polytope.
(See Kirwan \cite{kirwan;convexity-III}.)  Observe also that $\ca O$
is a symplectic quotient of the cotangent bundle $T^*K$.  This implies
that $\Delta(\ca O)$ is equal to the intersection of $\Delta(T^* K)$
with an affine subspace, where $\Delta(T^*K)$ denotes the moment cone
of $T^*K$ with respect to the $\ti K\times K$-action.  We shall in
fact write a complete set of inequalities for the cone $\Delta(T^*K)$.

We do so in terms of the Schubert bases of the cohomology groups
$H^\bu(X)$ and $H^\bu(\ti X)$, the canonical homomorphism $\phi\colon
H^\bu(X)\rightarrow H^\bu(\ti X)$, and the action on $H^\bu(X)$ of a
certain subset of the Weyl group of $K$, which we name the
\emph{relative Weyl set}.  Here $X$ and $\ti X$ denote the flag
varieties of $K$ and $\ti K$, respectively.  Determining the relative
Weyl set and the matrix of $\phi$ relative to the Schubert bases are
interesting combinatorial problems, which we do not know how to solve
in general.  In many examples, however, they can be solved and we
obtain explicit inequalities for the moment cone.  We treat a few of
these examples in detail.

An intriguing multiplicative version of Klyachko's theorem was
obtained by Agnihotri and Woodward
\cite{agnihotri-woodward;eigenvalues} and Belkale
\cite{belkale;local-systems}, who considered conjugacy classes in
$\U(n)$ and found eigenvalue inequalities for products of unitary
matrices in terms of quantum Schubert calculus.  An interesting
question is to what extent their results can be generalized to the
setting of the present paper.

Heckman showed that the geometric problem \eqref{item;geometric} is
closely related to the following algebraic problem.

\begin{enumerate}
\addtocounter{enumi}{1}
\item\label{item;algebraic}
Let $V$ be an irreducible $K$-module.  As a $\ti K$-module it breaks
up into isotypical components.  Which irreducible $\ti K$-modules
occur in $V$?
\end{enumerate}

Indeed, assume that $\ca O$ is \emph{integral} in the sense that it is
the orbit through an integral point $\lambda\in\t^*_+$ and that
$V_\lambda$ is the irreducible $K$-module with highest weight
$\lambda$.  Likewise, let $\ti V_\tlambda$ be the irreducible $\ti
K$-module attached to an integral coadjoint $\ti K$-orbit $\ti\ca O$.
Then if $\ti V_\tlambda$ occurs in $V_\lambda$, $\ti\ca O$ must be
contained in $f^*(\ca O)$.  (The converse question---which integral
$\ti K$-orbits inside $f^*(\ca O)$ correspond to irreducible
components of $V_\lambda$?---is in general much harder, and will be
discussed only briefly in this paper.)  As was pointed out by
Guillemin and Sternberg \cite{guillemin-sternberg;convexity;;1982},
Heckman's work implies that the following asymptotic version of
\eqref{item;algebraic} is equivalent to \eqref{item;geometric}.

\begin{enumerate}
\renewcommand\theenumi{$\text{ii}'$}
\item\label{item;asymptotic}
For which \emph{rational} points $(\tlambda,\lambda)$ in
$\ti\t^*_+\times\t^*_+$ does there exist a positive integer $n$ such
that $(n\tlambda,n\lambda)$ is integral and the irreducible module
$\ti V_{n\tlambda}$ occurs in $V_{n\lambda}$?
\end{enumerate}

Let $\ca S$ denote the set of all pairs of dominant weights
$(\tlambda,\lambda)$ such that $\ti V_\tlambda$ occurs in $V_\lambda$,
and let $\ca S'$ denote the set of all pairs $(\tlambda,\lambda)$ in
$\ti\t^*_+\times\t^*_+$ which satisfy \eqref{item;asymptotic}.
Kr\"amer proved in \cite{kramer;untergruppen-linearen} that $\ca S$ is
closed under addition.  Therefore $\ca S'$ is equal to the convex hull
over the rationals of $\ca S$.  The correspondence between problems
\eqref{item;geometric} and \eqref{item;asymptotic} is then as follows:
$\ca S'$ is equal to the set of all rational points inside
$\Delta(T^*K)$, and $\Delta(T^*K)$ is equal to the closure of $\ca
S'$.  (See e.g.\ \cite[Theorem 7.6]{sjamaar;convexity}.)  This
correspondence plays an important role in the proof of our main
result.

In the remainder of this paper we will, for convenience, allow $f$ to
be a homomorphism with finite kernel.  This does not change the nature
of problems \eqref{item;geometric} and \eqref{item;asymptotic} in any
way, because the sets $\ca S$ and $\Delta(T^*K)$ depend only on the
Lie algebras of $K$ and $\ti K$.

Section \ref{section;sub} is a collection of prerequisite results on
subgroups and Weyl chambers.  The main result of this paper, Theorem
\ref{theorem;cone}, is stated in Section \ref{section;schubert} and
proved in Section \ref{section;semistable}.  Section
\ref{section;examples} contains the examples.  In Appendix
\ref{section;flag} we review some material on flag varieties and
Appendix \ref{section;notation} is a compendium of our notational
conventions.

We are grateful to Friedrich Knop for a number of helpful comments.

\section{Subgroups and Weyl chambers}\label{section;sub}

Throughout this paper $\ti K$ and $K$ denote compact connected Lie
groups and $f\colon\ti K\to K$ a Lie group homomorphism with finite
kernel.  We fix once and for all maximal tori $T$ of $K$ and $\ti T$
of $\ti K$ such that $\ti T=T\cap\ti K$.  This section covers a number
of auxiliary results which are needed in Section
\ref{section;schubert} to state our main result.  Some are elementary,
but are included for lack of a reference.  In Section
\ref{subsection;weyl} we show that the Weyl group of $\ti K$ can be
viewed as a subgroup of the Weyl group of $K$ in a number of different
ways.  In Section \ref{subsection;relative} we discuss how the
partition of the Cartan subalgebra $\t=\Lie T$ into Weyl chambers
induces a partition of the small Cartan $\ti\t=\Lie\ti T$ which
refines \emph{its} partition into Weyl chambers.  In Section
\ref{subsection;plethysm} we work out an example.

For the purpose of understanding Section \ref{section;schubert} the
most important parts of this section are
\ref{definition;compatible}--\ref{proposition;face}.  On a first
reading the reader may want to skim through these and through Example
\ref{example;su3}, and then move on to Section \ref{section;schubert}.

A comment on our notation: whenever $O$ is any object or structure
associated with the group $K$ in a natural way, the corresponding
object or structure associated with $\ti K$ is designated by $\ti O$.

\subsection{The Weyl groups}\label{subsection;weyl}

The main result of this section is Theorem \ref{theorem;lift}, which
says that there is a canonically defined family of injective
homomorphisms from $\tW$ into $W$, which is parametrized by a
certain subgroup $\barW$ of $W$.  Because the homomorphism $f$ has
finite kernel, $\ti K$ and $f(\ti K)$ have the same Weyl groups.  For
this reason, let us assume in this section that $f$ is injective and
identify $\ti K$ with $f(\ti K)$.

If $K$ acts on a set $M$ then the \emph{normalizer} of a subset $M'$
is the subgroup $\ca N_K(M')=\{\,k\in K\mid kM'=M'\,\}$.  We say that
a subset $K'$ of $K$ \emph{normalizes} $M'$ if $K'\subset\ca N_K(M')$.
The \emph{centralizer} of $M'$ is the subgroup $\ca
Z_K(M')=\bigcap_{m\in M'}K_m$, which is normal in $\ca N_K(M')$.
(Here $K_m$ denotes the stabilizer or isotropy subgroup of $m$.)  We
will call the quotient
$$
W_K(M') =\ca N_K(M')/\ca Z_K(M')
$$
the \emph{Weyl group} of $K$ relative to $M'$.  Taking $M=K$, acting
on itself by conjugation, we have a Weyl group $W_K(K')$ relative to
any subset $K'$ of $K$.  In particular, \emph{the} Weyl group of $K$
is $W=W_K(T)$.

\begin{lemma}\label{lemma;normal}
\begin{enumerate}
\item\label{item;normal}
Let $E$ and $F$ be closed subgroups of $K$ and let $E$ act on $K$ by
conjugation.  If $E$ normalizes $F$\upn, then it also normalizes the
commutator subgroup $[F,F]$\upn, the centre $C(F)$\upn, the
centralizer $\ca Z_K(F)$\upn, and the identity component $F^0$.
\item\label{item;multiply}
Let $T'$ be any closed subgroup of $T$ and let $L=\ca Z_K(T')$.  Then
\begin{gather}
\label{equation;bignormal}
\ca N_K(L) =L\cdot\bigl(\ca N_K(L)\cap\ca N_K(T)\bigr)
=L^0\cdot\bigl(\ca N_K(L)\cap\ca N_K(T)\bigr),\\
\label{equation;smallnormal}
\ca N_K(T') =L\cdot\bigl(\ca N_K(T')\cap\ca N_K(T)\bigr)
=L^0\cdot\bigl(\ca N_K(T')\cap\ca N_K(T)\bigr),
\end{gather}
and the inclusion $\ca N_K(T')\cap\ca N_K(T)\hookrightarrow\ca
N_K(T')$ induces an isomorphism
\begin{equation}\label{equation;subquotient}
\ca N_W(T')/\ca Z_W(T')\iso W_K(T').
\end{equation}
Hence $W_K(T')$ is finite and the order of $W_K(T')$ divides the order
of $W$.
\end{enumerate}
\end{lemma}

\begin{proof}
\eqref{item;normal} is trivial.

It is clear that $\ca N_K(L)$ contains $L\cdot\bigl(\ca N_K(L)\cap \ca
N_K(T)\bigr)$.  Conversely, consider any $g\in\ca N_K(L)$.  Then
$gTg\inv\subset L^0$ because $L$ contains $T$ and $T$ is connected.
But then $gTg\inv$ is a maximal torus of $L^0$, so there exists $x\in
L^0$ such that $gTg\inv=xTx\inv$, or equivalently, $n=x\inv g\in \ca
N_K(T)$.  But then also $n\in L^0\cdot\ca N_K(L)\subset\ca N_K(L)$.
We have shown that $g=xn$ with $x\in L^0$ and $n\in\ca N_K(L)\cap \ca
N_K(T)$, which proves \eqref{equation;bignormal}.  Observe that $\ca
N_K(T')\subset\ca N_K(L)$ by \eqref{item;normal}.  Hence, by
\eqref{equation;bignormal}, $\ca N_K(T')$ is contained in the
intersection of $L^0\cdot\bigl(\ca N_K(L)\cap \ca N_K(T)\bigr)$ and
$\ca N_K(T')$, which is equal to $L^0\cdot\bigl(\ca N_K(T')\cap\ca
N_K(T)\bigr)$.  The inclusion $L\cdot\bigl(\ca N_K(T')\cap\ca
N_K(T)\bigr)\subset\ca N_K(T')$ is obvious.  This proves
\eqref{equation;smallnormal}.  The isomorphism
\eqref{equation;subquotient} is Exercice 4 in Bourbaki \cite[Ch.\ 9,
\S2]{bourbaki;groupes-algebres}.  Let us give a hint.  It is easy to
show that
\begin{gather}
\ca N_W(T') =\bigl(\ca N_K(T')\cap\ca N_K(T)\bigr)\big/T,\notag\\
\label{equation;smallerweyl}
\ca Z_W(T') =W_L(T) =\ca N_L(T)/T =\bigl(L\cap\ca N_K(T)\bigr)\big/T.
\end{gather}
This implies that $\ca N_W(T')/\ca Z_W(T')\cong\bigl(\ca
N_K(T')\cap\ca N_K(T)\bigr)\big/\bigl(L\cap\ca N_K(T)\bigr)$.  The
latter group clearly injects into $W_K(T')$ and it follows from
\eqref{equation;smallnormal} that this injection is surjective.
\end{proof}

We can apply \eqref{equation;subquotient} to the maximal torus $\ti T$
of the subgroup $\ti K$.  Since $\tW$ is a subgroup of $W_K(\ti T)$
in a natural way, the upshot is that $\tW$ is a subquotient of $W$.
But more is true.  Let us denote the subgroup $\ca Z_K(\ti T)$ by
$\bar K$.  This is a connected subgroup that contains the maximal
torus $T$.  Its root system is equal to the subsystem $\bar R$ of all
roots $\alpha\in R$ which vanish on $\ti\t$.  Alternatively, $\bar K$
can be described as the centralizer of $\bar\s$: $\bar K=\ca
Z_K(\bar\s)=K_{\bar\s}$, where $\bar\s$ is the smallest (closed) face
of the Weyl chamber which contains $\ti\t\cap\t_+$.  In other words,
$\bar\s$ is the face of $\t_+$ which is perpendicular to $\bar R$,
$$
\bar\s =\bigcap_{\alpha\in\bar R}\bigl\{\,\xi\in\t_+\bigm|
\alpha(\xi)=0\,\bigr\}.
$$  
Let $C=C(K)$, $\ti C=C(\ti K)$ and $\bar C=C(\bar K)$ be the centres
of $K$, $\ti K$ and $\bar K$, respectively.  Observe that
$C\subset\bar C$ and $\ti C\subset\bar C^0$.  The identity component
of $\bar C$ is the torus whose Lie algebra $\bar\lie c$ is equal to
the linear span of $\bar\s$.  For the normalizer and centralizer of
$\bar C$ we have the following result.

\begin{lemma}\label{lemma;subnormal} 
\begin{gather}
\label{equation;centralizer1}
\ca Z_K(\bar C)=\ca Z_K(\bar C^0)=\bar K,\\
\label{equation;normalizer1}
\ca N_K(\bar C)=\ca N_K(\bar C^0)=\ca N_K(\bar K),\\
\label{equation;normalizer2}
\ca N_K(\ti T)\subset\ca N_K(\bar C),\\
\label{equation;normalizer3}
\ca N_{\ti K}(\bar C)=\ca N_{\ti K}(\bar C^0)=\ca N_{\ti K}(\ti T).
\end{gather}
\end{lemma}

\begin{proof}
We have
\begin{equation}\label{equation;centralizer2}
\ca Z_K(\bar C)\subset\ca Z_K(\bar C^0) =\ca Z_K(\bar\lie c)
\subset\ca Z_K(\ti\t) =\ca Z_K(\ti T)=\bar K.
\end{equation}
On the other hand, $\bar C$ is the centre of $\bar K$, so $\bar
K\subset\ca Z_K(\bar C)$ and therefore all the inclusions in
\eqref{equation;centralizer2} are equalities.  This proves
\eqref{equation;centralizer1}.  By Lemma
\ref{lemma;normal}\eqref{item;normal} we have $\ca N_K(\bar
K)\subset\ca N_K(\bar C)$, because $\bar C$ is the centre of $\bar K$.
Likewise, $\ca N_K(\bar C^0)\subset\ca N_K(\bar K)$, because $\bar K$
is the centralizer of $\bar C^0$, and also $\ca N_K(\bar C)\subset\ca
N_K(\bar C^0)$.  In sum, we have shown
$$
\ca N_K(\bar C)\subset\ca N_K(\bar C^0)\subset\ca N_K(\bar
K)\subset\ca N_K(\bar C),
$$
which proves \eqref{equation;normalizer1}.  The inclusion
\eqref{equation;normalizer2} follows from \eqref{equation;normalizer1}
plus the fact that $\ca N_K(\ti T)$ is contained in $\ca N_K\bigl(\ca
Z_K(\ti T)\bigr)=\ca N_K(\bar K)$.  The inclusion
\eqref{equation;normalizer3} follows from \eqref{equation;normalizer1}
plus the fact that $\ca N_K(\bar K)\cap\ti K=\ca N_{\ti K}(\ti T)$.
\end{proof}

This implies the following relationships among the various Weyl
groups.  Here $\barW=W_{\bar K}(T)$ denotes the Weyl group of $\bar
K$.

\begin{lemma}\label{lemma;inclusion}
\begin{enumerate}
\item\label{item;inclusion}
The inclusions $\ca N_{\ti K}(\ti T)\subset\ca N_K(\ti T)$ and
\eqref{equation;normalizer2} induce injective homomorphisms
\begin{equation}\label{equation;injective}
W\inj W_K(\ti T)\inj W_K(\bar C^0) =W_K(\bar C).
\end{equation}
\item
Let $\bar S$ be a base \upn(set of simple roots\upn) of the root
system $\bar R$.  Then
\begin{equation}\label{equation;semidirect}
\ca N_W(\bar R)\cong\ca N_W(\bar S)\ltimes\barW,
\end{equation}
and the projection $\ca N_W(\bar R)\to\ca N_W(\bar S)$ induces an
isomorphism
\begin{equation}\label{equation;iso}
j_{\bar S}\colon W_K(\bar C)\iso\ca N_W(\bar S).
\end{equation}
\end{enumerate}
\end{lemma}

\begin{proof}
The kernel of the map $\ca N_{\ti K}(\ti T)\to W_K(\ti T)$ is $\ca
Z_K(\ti T)\cap\ti K=\ca Z_{\ti K}(\ti T)=\ti T$, so the induced
homomorphism $W\to W_K(\ti T)$ is injective.  By
\eqref{equation;centralizer1}, the kernel of the map $\ca N_K(\ti
T)\to W_K(\bar C)$ is $\ca Z_K(\ti T)$, so the induced homomorphism
$W_K(\ti T)\to W_K(\bar C)$ is injective.  This proves
\eqref{item;inclusion}.

Observe that $\barW$ is normal in $\ca N_W(\bar R)$ and that $\ca
N_W(\bar S)\cap\barW=\{1\}$.  If $w\in\ca N_W(\bar R)$, then $w\bar S$
is a base of $\bar R$, so $\bar ww\bar S=\bar S$ for some $\bar
w\in\barW$, which implies $w\in\ca N_W(\bar S)\cdot\barW$.  This
proves \eqref{equation;semidirect}.  Inverting the isomorphism
\eqref{equation;subquotient} (with $T'=\bar C$) and using
\eqref{equation;smallerweyl} yields an isomorphism
\begin{equation}\label{equation;roottorus}
W_K(\bar C)\iso\ca N_W(\bar C)/\ca Z_W(\bar C)=\ca N_W(\bar C)/\barW.
\end{equation}
We assert that $\ca N_W(\bar C)=\ca N_W(\bar R)$.  This follows from
the observation that $\ca N_W(\bar C)=\ca N_W(\bar C^0)$ is the set of
Weyl group elements which preserve $\bar\lie c$, the intersection of
the root hyperplanes defined by $\bar R$.  Hence we obtain from
\eqref{equation;roottorus} an isomorphism $W_K(\bar C)\cong\ca
N_W(\bar R)/\barW$.  By \eqref{equation;semidirect} there is a
canonical projection $\ca N_W(\bar R)\to\ca N_W(\bar S)$, which
induces an isomorphism $\ca N_W(\bar R)/\barW\cong\ca N_W(\bar S)$.
This proves \eqref{equation;iso}.
\end{proof}

The base $\bar S$ is unique up to an element of $\barW$, and a
different choice of base conjugates the isomorphism
\eqref{equation;iso}: $j_{\bar w\bar S}(w)=\bar wj_{\bar S}(w)\bar
w\inv$ for $\bar w\in\barW$ and $w\in W_K(\bar C)$.  We can summarize
these results as follows.

\begin{theorem}\label{theorem;lift}
A choice of a base $\bar S$ of the root subsystem $\bar R$ gives rise
to a diagram of canonically defined homomorphisms
\begin{equation}\label{equation;inclusion}
\tW\inj W_K(\ti T)\inj W_K(\bar C)\overset{j_{\bar S}}{\iso}\ca
N_W(\bar S) \inj\ca N_W(\bar C)\inj W.
\end{equation}
Hence the action of $\tW$ on $\ti T$ extends to an action on the torus
$T$ which preserves the intermediate subgroup $\bar C$ and the base
$\bar S$.  A different choice of base changes the lifting homomorphism
$j_{\bar S}$ by conjugation with an element of $\barW =\ca Z_W(\bar
C)$.  Therefore the extension of the $\tW$-action from $\ti T$ to
$\bar C$ is unique\upn; the extension from $\bar C$ to $T$ is unique
up to an element of $\barW$.
\qed
\end{theorem}

Via the isomorphism $j_{\bar S}$ we can regard $\tW$ as a subgroup of
$W$.  The image $j_{\bar S}(\tW)$ depends on $\bar S$, but the
subgroup $j_{\bar S}(\tW)\barW=\barW j_{\bar S}(\tW)$ does not.  We
usually abbreviate $j_{\bar S}$ to $j$ when the base $\bar S$ is
understood.

\begin{corollary}\label{corollary;lift}
The inclusions $\ti\t\to\bar\lie c\to\t$ and the dual projections
$\t^*\to\bar\lie c^*\to\ti\t^*$ are $\tW$-equivariant \upn(with
respect to any choice of $\bar S$\upn) and $\barW$-invariant.  The set
of projected roots $f^*(R)\subset\ti\t^*$ is $\tW$-stable.  The
embedding $\ti K/\ti T\to K/T$ induced by $f$ is $\tW$-equivariant.
\qed
\end{corollary}

\begin{example}\label{example;regular}
Observe that $\barW=\{1\}$ if and only if $\ti K$ contains a regular
element of $K$.  This is for instance the case if $\ti K$ is the unit
component of the fixed-point group of an automorphism of $K$.  (See
Bourbaki \cite[Ch.\ 9, \S5.3]{bourbaki;groupes-algebres}.)  If $\barW
=\{1\}$, then $\bar\s=\t_+$, $\bar K=T$, and $\bar C=\bar C^0=T$, so
\eqref{equation;inclusion} amounts to
$$
\tW\subset W_K(\ti T)\subset W_K(T)=\ca N_W(\emptyset)=W,
$$
where the inclusions are canonical.
\end{example}

\subsection{The relative Weyl set}\label{subsection;relative}

The Weyl chambers of the Cartan subalgebra $\t$ form a simplicial
subdivision in the sense that they cover the whole of $\t$, every
chamber is a simplicial cone, and the intersection of any two chambers
is a face of each.  In this section we study how the subdivision of
$\t$ into chambers induces a subdivision of $\ti\t$ which is
$\tW$-invariant and refines the chamber subdivision of $\ti\t$.

By a \emph{face} of $\t$ we mean a closed face of any Weyl chamber in
$\t$.  If $\s$ is a face, let $R_\s$ be the set of all $\alpha\in R$
such that $\alpha\ge0$ on $\s$.  Let $\g$ and $\h$ be the
complexifications of $\k$, resp.\ $\t$.  Then the subalgebra $\lie
p_\s =\h\oplus\bigoplus_{\alpha\in R_\s}\g_\alpha$ is parabolic, and
the correspondence $\s\mapsto\lie p_\s$ is a bijection between faces
of $\t$ and parabolic subalgebras of $\g$ containing $\h$.  We denote
the relative interior of a polyhedral subset $Q$ of a real vector
space by $Q^\circ$.

\begin{lemma}\label{lemma;face}
\begin{enumerate}
\item\label{item;subroot}
$\ti R\subset f^*(R)$.
\end{enumerate}
Let $\s$ and $\ti\s$ be faces of $\t$ and $\ti\t$\upn, respectively.
\begin{enumerate}
\addtocounter{enumi}{1}
\item\label{item;face}
$\s\cap\ti\t\subset\ti\s$ if and only if
$(\s\cap\ti\t)^\circ\cap\ti\s$ is nonempty.
\item\label{item;weylinclude}
If $\s\cap\ti\t\subset\ti\s$ then $w\s\cap\ti\t=\s\cap\ti\t$ for all
$w\in j(\tW_{\ti\s})\barW$\upn, where $\barW$ is the Weyl group of
$\bar K=\ca Z_K\bigl(f(\ti T)\bigr)$\upn, and $j$ is any lifting
homomorphism as in Theorem {\rm\ref{theorem;lift}}.
\end{enumerate}
\end{lemma}

\begin{proof}
\eqref{item;subroot} follows immediately from the fact that the
root-space decomposition $\g=\h\oplus\bigoplus_{\alpha\in R}\g_\alpha$
is a refinement of the weight-space decomposition of $\g$ relative to
the subtorus $\ti T$:
$$
\g =\bar\g\oplus\bigoplus_{\tlambda\in\ti\Lambda^*\setminus\{0\}}
\g_\tlambda,
$$
where $\ti\Lambda^*$ denotes the weight lattice in $\ti\t^*$ and
$\g_\tlambda=\bigoplus_{f^*(\alpha)=\tlambda}\g_\alpha$.

Let $\s$ be any face of $\t$.  The subset
$\{\,\alpha\mid\text{$\alpha\ge0$ on $\s\cap\ti\t$}\,\}$ of $R$
contains $R_\s$ and is therefore of the form $R_{\lie r}$, where $\lie
r$ is a face of $\s$.  It is clear that $\lie r$ is the smallest face
of $\t$ which contains $\s\cap\ti\t$, and therefore
$(\s\cap\ti\t)^\circ\subset\lie r^\circ$.

Now assume that $(\s\cap\ti\t)^\circ\cap\ti\s$ is nonempty.  Pick
$\tlambda\in(\s\cap\ti\t)^\circ\cap\ti\s$.  Then $\tlambda\in\lie
r^\circ$, so a root $\alpha$ is in $R_{\lie r}$ if and only if
$\alpha(\tlambda)\ge0$.  Furthermore $\tlambda\in\ti\s$, so if
$\ti\alpha\in\ti R_{\ti\s}$ then $\ti\alpha(\tlambda)\ge0$.  According
to \eqref{item;subroot} we can write $\ti\alpha=f^*(\alpha)$ with
$\alpha\in R$.  Hence $\alpha(\tlambda) =\ti\alpha(\tlambda)\ge0$,
which implies $\alpha\in R_{\lie r}$.  Therefore, if $\tmu$ is an
arbitrary element of $\s\cap\ti\t$, then
$\ti\alpha(\tmu)=\alpha(\tmu)\ge0$.  This shows that $\s\cap\ti\t$ is
contained in $\ti\s$.  Conversely, if $\s\cap\ti\t\subset\ti\s$ it is
obvious that $(\s\cap\ti\t)^\circ\cap\ti\s$ is nonempty.  This proves
\eqref{item;face}.

\eqref{item;weylinclude} follows from the fact that $\barW$ acts
trivially on $\ti\t$ and that $f_*$ is $\tW$-equivariant.
\end{proof}

Select a positive Weyl chamber $\ti\t_+$ in $\ti\t$.  Among all
chambers of $\t$, those which intersect $\ti\t_+$ in a cone of maximal
dimension are of special importance.

\begin{definition}\label{definition;compatible}
A Weyl chamber $\t_+$ in $\t$ is \emph{compatible} with $\ti\t_+$ if
$\dim\ti\t_+\cap\t_+=\dim\ti\t_+$.
\end{definition}

Equivalently, $\t_+$ is compatible with $\ti\t_+$ if there exists a
$\ti\xi_0$ in the interior of $\ti\t_+$ such that $\alpha(\ti\xi)\ge0$
for all $\alpha\in S$ and for all $\ti\xi\in\ti\t$ sufficiently close
to $\ti\xi_0$.  It is obvious that compatible chambers exist: pick an
arbitrary set of positive roots $\bar R_+$ in $\bar R$, pick
$\ti\xi_0\in\ti\t_+^\circ$, and let $R_+'$ be the set of $\alpha\in R$
such that $\alpha(\ti\xi_0)>0$.  Then $R_+=\bar R_+\cup R_+'$ is a set
of positive roots, and the chamber in $\t$ which is positive with
respect to $R_+$ is a compatible chamber.  Moreover, every compatible
chamber arises in this way.

Henceforth we shall fix chambers $\ti\t_+$ and $\t_+$ such that $\t_+$
is compatible with $\ti\t_+$.  We denote the corresponding positive
Borel subgroups of $\ti G=\ti K^\C$ and $G=K^\C$ by $\ti B$ and $B$,
respectively.

\begin{remark}\label{remark;choice}
This choice entails a choice of bases of the root systems $\ti R$ and
$R$, and hence also of a base of $\bar R$ and a lifting homomorphism
$j$ as in Theorem \ref{theorem;lift}.
\end{remark}

For arbitrary $w\in W$ we define $\ti\t_w$ to be the cone
$w\t_+\cap\ti\t$.  The following result gives a set of spanning
vectors for the cones dual to $\ti\t_w$ (and hence inequalities for
the $\ti\t_w$).  Let $\ca C$ denote the cone spanned by the positive
roots $R_+$.  (See Appendix \ref{section;notation} for our conventions
regarding roots and weights.)

\begin{lemma}\label{lemma;dualcone}
For all $w\in W$ the cone $\ti\t_w$ is dual to the cone $f^*(w\ca C)$.
\end{lemma}

\begin{proof}
Let $\ti\xi\in\ti\t$.  Then $f^*(w\alpha)(\ti\xi)\ge0$ for all
$\alpha\in R_+$ if and only if $\alpha(w\inv\ti\xi)\ge0$ for all
$\alpha\in R_+$.  This is equivalent to $w\inv\ti\xi\in\t_+$, i.e.\
$\ti\xi\in w\t_+\cap\ti\t=\ti\t_w$.
\end{proof}

\begin{definition}\label{definition;wcom}
The \emph{compatible Weyl set} $W\com$ is the set of all $w\in W$ such
that $w\t_+$ is compatible with $\ti\t_+$.
\end{definition}

Clearly $1\in W\com$ because $\t_+$ is assumed to be compatible with
$\ti\t_+$, and $\barW W\com=W\com$ because of Lemma
\ref{lemma;face}\eqref{item;weylinclude}.  Here are some further
helpful properties.

\begin{proposition}\label{proposition;plusface}
Let $w\in W\com$.  Then
\begin{enumerate}
\item\label{item;barw}
$\ti\t_{\bar ww}=\ti\t_w$ and $f^*(\bar ww\ca C)=f^*(w\ca C)$ for
$\bar w\in\barW$\upn,
\item\label{item;cone}
$\ti\t_w\subset\ti\t_+$ and $\ti\ca C\subset f^*(w\ca C)$\upn,
\item\label{item;uproot}
$(f^*)\inv(\ti R_+)\cap R\subset wR_+$ and $(f^*)\inv(\ti R_-)\cap
R\subset wR_-$\upn,
\item\label{item;downroot}
$\ti R_+\subset f^*(wR_+)$ and $\ti R_-\cap f^*(wR_+)=\emptyset$\upn,
\item\label{item;borel}
$\ti B=wBw\inv\cap\ti G$\upn,
\item\label{item;proper}
$f^*(w\ca C)$ is a proper cone\upn: $f^*(w\ca C)\cap-f^*(w\ca C)
=\{0\}$.
\end{enumerate}
\end{proposition}
 
\begin{proof}
\eqref{item;barw} follows from Lemma
\ref{lemma;face}\eqref{item;weylinclude}.

The first statement in \eqref{item;cone} is proved by setting
$\s=w\t_+$ and $\ti\s=\ti\t_+$ in Lemma
\ref{lemma;face}\eqref{item;face}.  The second statement then follows
from Lemma \ref{lemma;dualcone}.

Consider $\ti\alpha\in\ti R_+$ and $\alpha\in R$ such that
$f^*(\alpha)=\ti\alpha$.  Pick $\ti\xi\in w\t_+\cap\ti\t_+^\circ$;
then $\ti\xi$ can be written as $w\xi$ with $\xi\in\t_+$.  Hence
\begin{equation}\label{equation;maproot}
(w\inv\alpha)(\xi)=\alpha(w\xi)=f^*(\alpha)(\ti\xi)
=\ti\alpha(\ti\xi)>0,
\end{equation}
and therefore $w\inv\alpha\in R_+$, that is $\alpha\in wR_+$.  This
proves that $(f^*)\inv(\ti R_+)\cap R\subset wR_+$.  If $\ti\alpha$ is
negative, we get $(w\inv\alpha)(\xi)<0$ in \eqref{equation;maproot},
so we see that $(f^*)\inv(\ti R_-)\cap R\subset wR_-$.  This proves
\eqref{item;uproot}.

\eqref{item;downroot} follows immediately from \eqref{item;uproot} and
Lemma \ref{lemma;face}\eqref{item;subroot}.

The first assertion in \eqref{item;downroot} implies that $\ti\lie
b\subset w\lie b$ and therefore $\ti B\subset wBw\inv\cap\ti G$.  This
also shows that $wBw\inv\cap\ti G$ is a parabolic, and hence
connected, subgroup of $\ti G$, so to prove the reverse inclusion we
need only show that $w\lie b\cap\ti\g\subset\ti\lie b$.  This amounts
to showing that the parabolic subalgebra $w\lie b\cap\ti\g$ contains
no negative root spaces, which follows from the second statement in
\eqref{item;downroot}.  This proves \eqref{item;borel}.

Because of Lemma \ref{lemma;dualcone} the cone $f^*(w\ca C)$ is proper
if and only if $\dim\ti\t_w=\dim\ti\t$.  The latter is true because
$\ti\t_w=w\t_+\cap\ti\t_+$ by \eqref{item;cone} and $\dim
w\t_+\cap\ti\t_+=\dim\ti\t$ by assumption.  This proves
\eqref{item;proper}.
\end{proof}

Statement \eqref{item;proper} obviously fails for any $w\in W$ such
that $\ti\t_w$ is not top-dimensional.  On the other hand,
\eqref{item;cone}--\eqref{item;borel} are true, with exactly the same
proof, under the more general hypothesis that
$\ti\t_w^\circ\cap\ti\t_+^\circ$ be nonempty.  Therefore
\eqref{item;cone}--\eqref{item;borel} are necessary, but not
sufficient, conditions for $w$ to be a compatible element.  The proof
also shows that \eqref{item;cone}--\eqref{item;borel} fail when
$\ti\t_w$ is not contained in $\ti\t_+$.  (When $\ti\t_w$ is contained
in the \emph{boundary} of $\ti\t_+$ then \eqref{item;cone} is true but
\eqref{item;uproot}--\eqref{item;borel} may fail, depending on $w$.)

According to \eqref{item;barw} the top-dimensional cones $\t_w$ are
labelled by the left quotient $\barW\backslash W\com$.  It is
well-known that every coset $\barW w$ contains a unique element of
minimal length, called its \emph{shortest representative}.  A more
convenient set of labels is obtained by selecting from each coset the
shortest representative.

\begin{definition}\label{definition;relative}
The \emph{relative Weyl set} $W\rel$ is the set of shortest
representatives of the left quotient $\barW\backslash W\com$.  For
$v\in W\rel$ the cone $\ti\t_v$ is called a \emph{cubicle} in
$\ti\t_+$.
\end{definition}

The sets $W\com$ and $W\rel$ depend of course on the choice of $\t_+$
and $\ti\t_+$.  We emphasize that in general neither $W\com$ nor
$W\rel$ are subgroups of $W$.  We now show that the collection of
cubicles forms a well-behaved conical subdivision of $\ti\t_+$ which
is in one-to-one correspondence with $W\rel$.

\begin{proposition}\label{proposition;face}
\begin{enumerate}
\item\label{item;subdivision}
The chamber $\ti\t_+$ is the union of the cubicles\upn,
\begin{equation}\label{equation;subdivision}
\ti\t_+=\bigcup_{v\in W\rel}\ti\t_v,
\end{equation}
and the intersection of any two cubicles is a face of each.
\item\label{item;overlap}
If $\ti\t_{v_1}^\circ\cap\ti\t_{v_2}\ne\emptyset$ for some $v_1$ and
$v_2\in W\rel$\upn, then $v_1=v_2$.
\item\label{item;wrel}
$W\rel=\{1\}$ if and only if $\ti\ca C=f^*(\ca C)$.
\end{enumerate}
\end{proposition}

\begin{proof}
The inclusion $\bigcup_{v\in W\rel}\ti\t_v\subset\ti\t_+$ follows from
Proposition \ref{proposition;plusface}\eqref{item;cone}.  Because $\t$
is the union of its chambers $w\t_+$, the positive chamber $\ti\t_+$
of $\ti\t$ is the union of all cones $\ti\t_w$ with
$\ti\t_w^\circ\cap\ti\t_+\ne\emptyset$.  From this finite union we can
clearly delete those cones whose dimension is less than $\dim\ti\t$,
which leaves only the cubicles.  The fact that the intersection of two
cubicles is a face of each follows from the corresponding fact for the
subdivision $\bigcup_{w\in W}w\t_+$ of $\t$.  This proves
\eqref{item;subdivision}.

If $\ti\t_{v_1}^\circ\cap\ti\t_{v_2}\ne\emptyset$, then
$\ti\t_{v_1}=\ti\t_{v_2}$ by \eqref{item;subdivision}.  The cone
$\ti\t_{v_1}=\ti\t_{v_2}$ has an open relative interior in $\ti\t_+$
and furthermore $W_{\ti\xi}=W_{\bar\lie c}=\barW$ for an open dense
set of $\ti\xi\in\ti\t_+^\circ$.  This implies that there exists
$\ti\xi$ in $v_1\t_+\cap\ti\t_+^\circ =v_2\t_+\cap\ti\t_+^\circ$ such
that $W_{\ti\xi} =\barW$.  This in turn implies that the two faces
$v_1\t_+\cap\bar\lie c$ and $v_2\t_+\cap\bar\lie c$ of $\t$ are equal
to one another.  Let $w=v_2v_1\inv$; then $w$ maps the chamber
$v_1\t_+$ to the chamber $v_2\t_+$ and therefore must fix their
intersection $v_1\t_+\cap\bar\lie c=v_2\t_+\cap\bar\lie c$.  But then
$w$ fixes $\bar\lie c$, so $w\in\barW$, i.e.\ $v_2\in v_1\barW$.  As
$v_1$ and $v_2$ are both shortest representatives of their cosets
$\bmod\barW$, we conclude that $v_1=v_2$.  This proves
\eqref{item;overlap}.

Dualizing the condition $\ti\ca C=f^*(\ca C)$ we obtain
$\ti\t_+=\ti\t_1$.  The latter condition is equivalent to the
subdivision \eqref{equation;subdivision} consisting of one cubicle
only.  According to \eqref{item;overlap} this means that
$W\rel=\{1\}$.  This proves \eqref{item;wrel}.
\end{proof}

We show in Section \ref{subsection;duality} that the subdivision
\eqref{equation;subdivision} is symmetric under the standard (duality)
involution of $\ti\t_+$.

Observe that the cubicles are proper cones if and only if $\ti T\cap
C$ is discrete.  Indeed, every cubicle $\ti\t_v$ contains the linear
subspace $\ti\t\cap\lie c=\Lie(\ti T\cap C)$.  On the other hand,
$\lie c$ is the largest subspace contained in $v\t_+$ for all $v$, so
if $\ti\t\cap\lie c=0$, then $\ti\t_v$ is a proper cone.

It is clear from Proposition \ref{proposition;face} that by
intersecting all the chambers in $\ti\t$ with chambers in $\t$ one
obtains a subdivision of the entire Cartan subalgebra $\ti\t$ into
polyhedral cones which refines the subdivision $\ti\t=\bigcup_{\ti
w\in\tW}\ti w\ti\t_+$ and which, by Corollary \ref{corollary;lift}, is
$\tW$-invariant.

An illustration of Propositions \ref{proposition;plusface} and
\ref{proposition;face} can be found in Section
\ref{subsection;plethysm}.  For the applications in Section
\ref{section;schubert} it would be highly desirable to find the
inequalities that describe each of the cones $f^*(v\ca C)$, or
equivalently, find an explicit set of spanning vectors for each of the
cubicles $\ti\t_v$.  This we have not been able to do in general.

We conclude this section with some miscellaneous observations.  Let
$W'$ be the set of all $w\in W$ such that $(f^*)\inv(\ti R_+)\cap
R\subset wR_+$.  By Proposition \ref{proposition;plusface} and the
remarks following it, $W'$ includes $W\com$ as well as all other $w$
such that $\ti\t_w$ contains interior points of $\ti\t_+$.  Put
$R'=\bigcap_{w\in W'}wR_+$ and $\ca C'=\bigcap_{w\in W'}w\ca C$.

\begin{proposition}\label{proposition;relroot}
\begin{enumerate}
\item\label{item;intersectroot}
$\ti R_+\subset f^*(R')$.
\item\label{item;intersect}
$\ti\ca C=f^*(\ca C')$.
\item\label{item;simple}
$\bar S=vS\cap\bar R$ for all $v\in W\rel$.
\end{enumerate}
\end{proposition}

\begin{proof}
\eqref{item;intersectroot} follows immediately from the definition of
$R'$.

For any subset $\eu S$ of a vector space $\eu E$ we define $\eu
S\spcheck$ to be the cone in $\eu E^*$ consisting of all functionals
which are nonnegative on $\eu S$.  Then $\eu S_1\spcheck\cap\eu
S_2\spcheck =(\eu S_1\cup\eu S_2)\spcheck$ for all subsets $\eu S_1$
and $\eu S_2$.  Applying this observation to
\eqref{equation;subdivision} we obtain $\ti\ca C=\bigcap_{v\in
W\rel}\ti\t_v\spcheck$.  By Lemma \ref{lemma;dualcone} this implies
$\ti\ca C=\bigcap_{v\in W\rel}f^*(v\ca C)$, and hence $f^*(\ca
C\rel)\subset\ti\ca C$.  The reverse inclusion follows from
\eqref{item;intersectroot}.  This proves \eqref{item;intersect}.

For $w\in W$ let $\bar\lie c_w$ be the cone $w\t_+\cap\bar\lie c$.  We
assert that for all $w\in W$
\begin{equation}\label{equation;base}
\text{$\dim\bar\lie c_w=\dim\bar\lie c\iff wS\cap\bar R$ is a base of
$\bar R$}.
\end{equation}
Indeed, if $\dim\bar\lie c_w=\dim\bar\lie c$, then $\bar\lie c$ is the
linear span of $\bar\lie c_w$.  Furthermore, $wR_+\cap\bar R$ is
clearly a set of positive roots for $\bar R$.  Let $\gamma\in
wR_+\cap\bar R$ and write $\gamma=w\beta$ with $\beta\in R_+$.  Let
$\alpha\in S$ be a simple root that contributes to $\beta$.  Then
$\alpha\ge0$ on $\t_+$, so $w\alpha\ge0$ on $\bar\lie c_w$.  On the
other hand, $\gamma=w\beta\in\bar R$, so $w\beta$ vanishes on
$\bar\lie c$.  This is only possible if $w\alpha$ vanishes on
$\bar\lie c_w$.  Hence $w\alpha$ vanishes on $\bar\lie c$ and
therefore $w\alpha\in\bar R$.  We conclude that $\gamma$ is a linear
combination with nonnegative coefficients of elements in $wS\cap\bar
R$, so this set forms a base of $\bar R$.  Conversely, if $wS\cap\bar
R$ is a base of $\bar R$, then $\bar\lie c$, resp.\ $\bar\lie c_w$, is
equal to the product $\lie c\times\lie d$, where $\lie c$ is the
centre of $\k$ and $\lie d$ is the linear, resp.\ positive, span of
all fundamental coweights that are perpendicular to $wS\cap\bar R$.
Hence $\dim\bar\lie c_w=\dim\bar\lie c$, which finishes the proof of
\eqref{equation;base}.

Now observe that $\dim\bar\lie c_v =\dim\bar\lie c$ because $v$ is
compatible.  We infer from \eqref{equation;base} that $vS\cap\bar R$
is a base of $\bar R$.  Therefore $vS\cap\bar R=\bar v\inv\bar S$ for
some $\bar v\in\barW$, or equivalently, $\bar wvS\cap\bar R=\bar S$.
This implies that
\begin{equation}\label{equation;positive}
\bar wvR_+\cap\bar R=\bar R_+.
\end{equation}
Take any $\alpha\in R_+$.  If $v\alpha\in\bar R$, then $\bar
wv\alpha\in\bar R_+$ because of \eqref{equation;positive}.  If
$v\alpha$ is not in $\bar R$, then $v\alpha\in R_+$ if and only if
$\bar wv\alpha\in R_+$, because $\bar w\in\barW$ permutes the elements
of $R_+\setminus\bar R_+$.  This shows that the length of $v$, which
is equal to the cardinality of $vR_+\cap R_-$, is greater than or
equal to the length of $\bar wv$.  But then $\bar w=1$, because $v$ is
the shortest element in $\barW v$.  The result is that $vS\cap\bar
R=\bar S$.
\end{proof}

The condition \eqref{item;simple} is not sufficient for $v$ to be in
$W\rel$ for the following reason.  The proof shows that
\eqref{item;simple} holds for any $v\in W$ which is the shortest
$\barW$-representative of a $w\in W$ such that $\dim\bar\lie c_w
=\dim\bar\lie c$.  Let us take $w$ such that $\ti\t_w$ is compatible
with $\ti w\ti\t_+$, where $\ti w$ is an arbitrary element of $\tW$.
Then $\dim\bar\lie c_w =\dim\bar\lie c$, but if $\ti w\ne1$ then
$w\not\in W\rel$, and therefore the shortest representative of $w$ is
not in $W\rel$.

\subsection{Plethysms}\label{subsection;plethysm}

Interesting examples are provided by representation theory.  A unitary
representation of an arbitrary compact connected group $\ti K$ on a
finite-dimensional Hermitian vector space $V$ is nothing but a
homomorphism $f\colon\ti K\to K$, where $K$ is the unitary group
$\U(V)$.  Note that every nontrivial representation has finite kernel
if $\ti K$ is simple.
We begin by stating some generalities concerning this class of
examples and then work out a simple case in more detail.

Let $V=\bigoplus_{\mu\in\ti\Lambda^*}V^\mu$ be the weight-space
decomposition of $V$ and let $n(\mu)$ denote the multiplicity of a
weight $\mu\in\ti\Lambda^*$.  Then $\sum n(\mu) =n=\dim V$.  The
centralizer $\bar K=\ca Z_K(\ti T)$ and its centre $\bar C$ are given
by
$$
\bar K=\prod_\mu\U(V^\mu),\qquad \bar C =\prod_\mu C^\mu,
$$
where $C^\mu$ denotes the one-dimensional centre of $\U(V^\mu)$.  Let
us choose a maximal torus $T$ of $\U(V)$ which contains $f(\ti T)$.
This boils down to a choice of (unordered) bases $\ca B^\mu
=\bigl\{e_1^\mu,e_2^\mu,\dots,e_{n(\mu)}^\mu\bigr\}$ of the weight
spaces $V^\mu$.  Let $\ca B$ be the basis
\begin{equation}\label{equation;basis}
\ca B=\coprod_\mu\ca B^\mu
\end{equation}
of $V$.  Then the Lie algebra $\t$ of $T$ can be identified with the
real span of $\ca B$ by sending $h\in\t$ to the vector
\begin{equation}\label{equation;vector}
\sum_{(\mu,k)}h(e_k^\mu)e_k^\mu.
\end{equation}
The Weyl group $W=S_n$ is the group of permutations of $\ca B$; $\ca
N_W(\bar C)$ consists of those permutations that preserve the
decomposition \eqref{equation;basis}; and $\barW=\ca Z_W(\bar C)
=\prod_\mu S_{n(\mu)}$ consists of those permutations that map each
$\ca B^\mu$ into itself.  Let us identify $\t$ with $\t^*$ by means of
the standard inner product on $\t$ determined by the basis $\ca B$.
Then the root system $R=R\spcheck$ is the set of all vectors of the
form $e_k^\mu-e_l^\nu$, where $(\mu,k)\ne(\nu,l)$, and the subsystem
$\bar R$ consists of the vectors $e_k^\mu-e_l^\mu$.

Observe that $\bar R\ne\emptyset$ and $\barW\ne\{1\}$ as soon as $V$
has a weight of multiplicity greater than $1$.  For any such
representation $V$ all points in $f(\ti K)$ are singular considered as
points in $\U(V)$.

The matrix of $f_*$ (relative to the coroots in $\ti\t$ and the basis
$\ca B$ in $\t$) and the set of projected roots $f^*(R)\subset\ti\t^*$
can be read off easily from the weight diagram of $V$.

\begin{lemma}\label{lemma;embedding}
Put $e^\mu=\sum_{i=1}^{n(\mu)}e_k^\mu$ for all weights $\mu$ of $V$.
Then the inclusion $f_*\colon\ti\t\hookrightarrow\bar\lie
c\hookrightarrow\t$ is given by
$$
f_*(\ti h)=\sum_\mu\mu(\ti h)e^\mu.
$$
Hence $f^*(e_k^\mu)=\mu$ and $f^*(e_k^\mu-e_l^\nu)=\mu-\nu$ for all
$(\mu,k)$ and $(\nu,l)$.
\end{lemma}

\begin{proof}
For all $(\mu,k)$ we have $f_*(\ti h)(e_k^\mu)=\mu(\ti h)e_k^\mu$
because $e_k^\mu$ has weight $\mu$.  Therefore, according to the
identification \eqref{equation;vector}, $f_*(\ti h)$ is equal to
$\sum_\mu\mu(\ti h)e^\mu$.

Dually, we have 
$$
\bigl\langle f^*(e_k^\mu),\ti h\bigr\rangle =\bigl\langle
e_k^\mu,f_*(\ti h)\bigr\rangle =\sum_\nu\nu(\ti h)\langle
e_k^\mu,e^\nu\rangle =\mu(\ti h),
$$
so $f^*(e_k^\mu)=\mu$.
\end{proof}

According to the discussion following Definition
\ref{definition;compatible}, a compatible pair of chambers $\ti\t_+$
and $\t_+$ is specified by a set of positive roots in $\bar R$, i.e.\
an ordering of each of the bases $\ca B^\mu$, and by a set of positive
roots in $\ti R$, which implies an ordering $\pre$ on the weight
lattice $\ti\Lambda^*$.  This gives rise to an ordering on the set of
pairs $(\mu,k)$ defined by $(\mu,k)\ge(\nu,l)$ if $\mu\suc\nu$ and
$(\mu,k)\ge(\mu,l)$ if $k\ge l$.  In terms of this ordering $R_+$ is
given by the set of all $e_k^\mu-e_l^\nu$ with $(\mu,k)\ge(\nu,l)$.

\begin{example}\label{example;su3}
Take $\ti K=\SU(3)$ and consider the unitary irreducible
representation $\ti V_\tlambda$ with highest weight
$\tlambda=2\ti\pi_1+\ti\pi_2$, where $\ti\pi_1$ and $\ti\pi_2$ are the
fundamental weights of $\ti K$.  Then $\dim\ti V_\tlambda=15$, so this
representation defines an embedding $f\colon\ti K\to K$, where
$K=\U(15)$.  We number the weights as in Figure \ref{figure;number}.
This ordering specifies a positive chamber $\t_+$ in $K$ which is
compatible with the (standard) positive chamber $\ti\t_+$ in $\ti K$.
We denote the corresponding bases of the root systems $R$ and $\ti R$
by $S$ and $\ti S$, respectively.

\begin{figure}
\setlength{\unitlength}{0.08mm}
$$
\begin{picture}(800,800)(-400,-400)
\thinlines
\dashline{10}(200,346.41016)(-200,-346.41016)
\dashline{10}(-400,0)(400,0)
\dashline{10}(-200,346.41016)(200,-346.41016)
\thicklines \path(0,0)(-150,86.60254) 
\path(-140,74)(-150,86.60254)(-134,85) 
\path(0,0)(150,86.60254) 
\path(140,74)(150,86.60254)(134,85) 
\put(-50,259.80762){\circle*{10}}
\put(100,173.20508){\circle*{10}}
\put(-200,173.20508){\circle*{10}}
\put(250,86.60254){\circle*{10}}
\put(-50,86.60254){\circle*{10}}
\put(100,0){\circle*{10}}
\put(-200,0){\circle*{10}}
\put(250,-86.60254){\circle*{10}}
\put(-50,-86.60254){\circle*{10}}
\put(100,-173.20508){\circle*{10}}
\put(-200,-173.20508){\circle*{10}}
\put(-50,-259.80762){\circle*{10}}
\put(-40,270){\makebox(0,0)[lb]{\smash{$\tlambda$}}}
\put(-100,90){\makebox(0,0)[lb]{\smash{$\ti\pi_1$}}}
\put(-170,50){\makebox(0,0)[lb]{\smash{$\ti\alpha_1$}}}
\put(160,65){\makebox(0,0)[lb]{\smash{$\ti\alpha_2$}}}
\put(-70,270){\makebox(0,0)[lb]{\smash{$\scriptstyle1$}}}
\put(-220,190){\makebox(0,0)[lb]{\smash{$\scriptstyle2$}}}
\put(80,190){\makebox(0,0)[lb]{\smash{$\scriptstyle3$}}}
\put(-50,100){\makebox(0,0)[lb]{\smash{$\scriptstyle4,5$}}}
\put(230,100){\makebox(0,0)[lb]{\smash{$\scriptstyle6$}}}
\put(-220,10){\makebox(0,0)[lb]{\smash{$\scriptstyle7$}}}
\put(100,10){\makebox(0,0)[lb]{\smash{$\scriptstyle8,9$}}}
\put(-50,-75){\makebox(0,0)[rb]{\smash{$\scriptstyle10,11$}}}
\put(220,-75){\makebox(0,0)[lb]{\smash{$\scriptstyle12$}}}
\put(-230,-160){\makebox(0,0)[lb]{\smash{$\scriptstyle13$}}}
\put(100,-155){\makebox(0,0)[lb]{\smash{$\scriptstyle14$}}}
\put(-80,-245){\makebox(0,0)[lb]{\smash{$\scriptstyle15$}}}
\end{picture}
$$
\caption{Numbering of weights for $\SU(3)$-module with highest weight
$\tlambda=2\ti\pi_1+\ti\pi_2$}
\label{figure;number}
\end{figure}
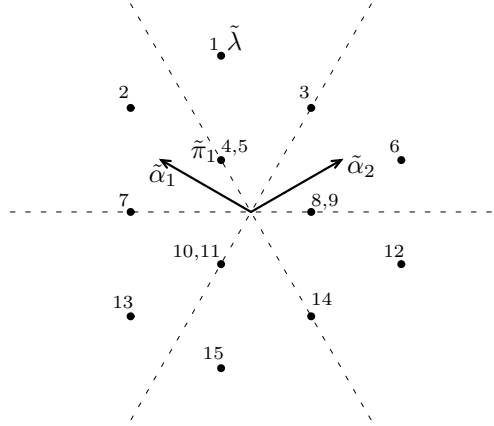

The projected roots are easy to compute from Lemma
\ref{lemma;embedding}, and by Lemma \ref{lemma;dualcone} the cubicles
in $\ti\t_+$ are then determined by the hyperplanes perpendicular to
the $f^*(\alpha)$.  The matrix of $f^*$ relative to the bases of
fundamental weights in $\SU(15)$ and $\SU(3)$ is
\begin{equation}\label{equation;matrix}
\setcounter{MaxMatrixCols}{14}
\begin{pmatrix}
2&5&5&6&7&5&7&6&5&5&5&2&3&1\\1&0&2&2&2&5&3&4&5&4&3&5&2&2
\end{pmatrix}
\end{equation}
The projected roots and the cubicles are shown in Figure
\ref{figure;peacock}, where we have used the trace form to identify
$\ti\t$ with its dual.  The arrows, with multiplicities, indicate the
images of the $105$ positive roots.  If an arrow occurs as the image
of one of the $14$ simple roots, the number of times it so occurs is
written in parentheses.  Three positive roots are mapped to $0$,
namely $\alpha_4$, $\alpha_8$ and $\alpha_{10}$, so that $\bar S=\bar
R_+ =\{\alpha_4,\alpha_8,\alpha_{10}\}$ and $\dim\bar C=14-3=11$.

\begin{figure}
\setlength{\unitlength}{0.1mm}
$$
\begin{picture}(800,800)(-400,-200)
\thinlines
\dashline{10}(0,0)(-100,-173.20508) \dashline{10}(-400,0)(400,0)
\dashline{10}(0,0)(100,-173.20508)
\thicklines
\texture{cccccccc 0 0 0 cccccccc 0 0 0 cccccccc 0 0 0 cccccccc 0 0 0
	cccccccc 0 0 0 cccccccc 0 0 0 cccccccc 0 0 0 cccccccc 0 0 0 }
\shade\path(0,0)(-100,519.61524)(-90,541)(-80,559)(-70,574)(-60,586)
(-50,595)(-40,601)(-30,604)(-20,604)(-10,602)(0,600)(0,0)
\thinlines
\texture{c0c0c0c0 0 0 0 0 0 0 0 c0c0c0c0 0 0 0 0 0 0 0 c0c0c0c0 0 0 0
	0 0 0 0 c0c0c0c0 0 0 0 0 0 0 0 }
\shade\path(0,0)(-300,519.61524)(-290,511)(-280,504)(-270,497)(-260,491)
(-250,485)(-240,480)(-230,475)(-220,471)(-210,467)(-200,464)(-190,462)
(-180,461)(-170,460)(-160,462)(-150,466)(-140,472)(-130,480)(-120,490)
(-110,501)(-100,519.61524)(0,0)
\shade\path(0,0)(0,600)(10,596)(20,590)(30,582)(40,572)(50,560)(60,547)
(70,536)(80,528)(90,523)(100,519.61524)(0,0)
\shade\path(0,0)(100,519.61524)(110,517)(120,515)(130,514)(140,513)
(150,513)(160,514)(170,516)(180,519)(190,523)(200,528)(210,533)(220,537)
(230,540)(240,542)(250,543)(260,544)(270,542)(280,537)(290,528)
(300,519.61524)(0,0)
\thicklines \path(0,0)(-150,86.60254) 
\path(-140,74)(-150,86.60254)(-134,85) 
\path(0,0)(150,86.60254) 
\path(140,74)(150,86.60254)(134,85) 
\path(0,0)(-300,173.20508)
\path(-290,161.6)(-300,173.20508)(-284,171.6)
\path(0,0)(300,173.20508) \path(290,161.6)(300,173.20508)(284,171.6)
\path(0,0)(-450,259.80762)
\path(-440,247.2)(-450,259.80762)(-434,258.2)
\path(0,0)(450,259.80762) \path(440,247.2)(450,259.80762)(434,258.2)
\path(0,0)(-450,86.60254) \path(-437,78)(-450,86.60254)(-434,90)
\path(0,0)(450,86.60254) \path(437,78)(450,86.60254)(434,90)
\path(0,0)(-300,346.41016) \path(-295,330)(-300,346.41016)(-285,342)
\path(0,0)(300,346.41016) \path(295,330)(300,346.41016)(285,342)
\path(0,0)(-150,433.0127) \path(-150,415)(-150,433.0127)(-140,420)
\path(0,0)(150,433.0127) \path(150,415)(150,433.0127)(140,420)
\path(0,0)(-150,259.80762) \path(-148,243)(-150,259.80762)(-137,250)
\path(0,0)(150,259.80762) \path(148,243)(150,259.80762)(137,250)
\path(0,0)(0,173.20508) \path(-5,160)(0,173.20508)(5,160)
\path(0,0)(0,346.41016) \path(-5,333)(0,346.41016)(5,333)
\path(0,0)(0,519.61524) \path(-5,504)(0,519.61524)(5,504)
\path(0,0)(300,0) \path(290,-5)(300,0)(290,5)
\put(0,0){\circle*{10}}
\put(-50,259.80762){\circle*{10}}
\put(-50,86.60254){\circle*{10}}
\put(50,86.60254){\circle*{10}}
\put(-40,270){\makebox(0,0)[lb]{\smash{$\tlambda$}}}
\put(-100,110){\makebox(0,0)[lb]{\smash{$\ti\pi_1$}}}
\put(75,70){\makebox(0,0)[lb]{\smash{$\ti\pi_2$}}}
\put(-170,55){\makebox(0,0)[lb]{\smash{$\ti\alpha_1$}}}
\put(160,65){\makebox(0,0)[lb]{\smash{$\ti\alpha_2$}}}
\put(-160,95){\makebox(0,0)[lb]{\smash{$\scriptstyle15(4)$}}}
\put(135,95){\makebox(0,0)[lb]{\smash{$\scriptstyle15$}}}
\put(-305,180){\makebox(0,0)[lb]{\smash{$\scriptstyle6$}}}
\put(295,182){\makebox(0,0)[lb]{\smash{$\scriptstyle6$}}}
\put(-460,266){\makebox(0,0)[lb]{\smash{$\scriptstyle1$}}}
\put(460,260){\makebox(0,0)[lb]{\smash{$\scriptstyle1$}}}
\put(-460,98){\makebox(0,0)[lb]{\smash{$\scriptstyle2(2)$}}}
\put(460,95){\makebox(0,0)[lb]{\smash{$\scriptstyle2$}}}
\put(-310,355){\makebox(0,0)[lb]{\smash{$\scriptstyle2$}}}
\put(310,345){\makebox(0,0)[lb]{\smash{$\scriptstyle2$}}}
\put(-170,425){\makebox(0,0)[lb]{\smash{$\scriptstyle2$}}}
\put(160,425){\makebox(0,0)[lb]{\smash{$\scriptstyle2$}}}
\put(-175,260){\makebox(0,0)[lb]{\smash{$\scriptstyle8$}}}
\put(160,260){\makebox(0,0)[lb]{\smash{$\scriptstyle8$}}}
\put(-30,175){\makebox(0,0)[lb]{\smash{$\scriptstyle15$}}}
\put(-20,347){\makebox(0,0)[lb]{\smash{$\scriptstyle6$}}}
\put(-20,520){\makebox(0,0)[lb]{\smash{$\scriptstyle1$}}}
\put(290,-28){\makebox(0,0)[lb]{\smash{$\scriptstyle8(5)$}}}
\put(-7,-28){\makebox(0,0)[lb]{\smash{$\scriptstyle3(3)$}}}
\put(-200,400){\makebox(0,0)[lb]{\smash{$\mathbf1$}}}
\put(-50,450){\makebox(0,0)[lb]{\smash{$\mathbf2$}}}
\put(30,450){\makebox(0,0)[lb]{\smash{$\mathbf3$}}}
\put(220,480){\makebox(0,0)[lb]{\smash{$\mathbf4$}}}
\end{picture}
$$
\caption{Peacock's tail with 105 feathers.  Projected positive roots
$f^*(R_+)$ and cubicles $\ti\t_v$ for embedding $\SU(3)\to\U(15)$
defined by dominant weight $\tlambda=2\ti\pi_1+\ti\pi_2$}
\label{figure;peacock}
\end{figure}

The centralizer $\bar K=\ca Z_K(\ti T)=\ca Z_K(\bar C)$ is a group of
block-diagonal matrices with nine $1\times1$-blocks and three
$2\times2$-blocks and hence is isomorphic to $\U(1)^9\times\U(2)^3$.
Its Weyl group is $\barW=(\Z/2\Z)^3$; it is generated by the simple
reflections $s_4$, $s_8$ and $s_{10}$, and acts by permuting rows and
columns of each of the $2\times2$-blocks.  The group $\ca N_W(\bar C)$
is isomorphic to $S_9\times(\Z/2\Z)^3\rtimes S_3$.  The copy of $S_9$
acts on $\bar K$ and its centre $\bar C$ by permuting the
$1\times1$-blocks, the copy of $S_3$ permutes the $2\times2$-blocks,
and the copy of $(\Z/2\Z)^3$ permutes the rows and columns of the
individual $2\times2$-blocks.  The group $\ca N_W(\bar S)$ is
isomorphic to $S_9\times S_3$.  It consists of all permutations $w$ in
$W\cong S_{15}$ that map the sets $\{1,2,3,6,7,12,13,14,15\}$ and
$\{4,8,10\}$ into themselves and satisfy $w(i)=i+1$ for $i=4$, $8$, or
$10$.

The homomorphism $j\colon\tW\to\ca N_W(\bar S)$ is found as follows.
By Theorem \ref{theorem;lift}, for each $\ti w\in\tW$ the image of
$\ti w$ under $j$ is the unique $w\in\ca N_W(\bar S)$ such that
$f_*(\ti w\ti\xi)=wf_*(\ti\xi)$ for all $\ti\xi\in\ti\t$.  In fact, in
this condition it suffices to take a single $\ti\xi$ in the interior
of one of the cones $\ti\t_v$, e.g.\ $\ti\xi
=6\ti\alpha\spcheck_1+5\ti\alpha\spcheck_2\in\ti\t_1$.  An easy
calculation now yields that on generators $j$ is given by
\begin{equation}\label{equation;j}
\begin{aligned}
j(\ti s_1) &=(1\;6)(2\;12)(4\;8)(5\;9)(7\;14)(13\;15),\\
j(\ti s_2) &=(1\;2)(3\;7)(6\;13)(8\;10)(9\;11)(12\;15).
\end{aligned}
\end{equation}
The other seven natural embeddings of $\tW$ into $W$ are obtained by
conjugating $j$ with any of the transpositions $(4\;5)$, $(8\;9)$,
$(10\;11)$, or a product thereof.

The relative Weyl set consists of the four elements
\begin{equation}\label{equation;relative}
\begin{gathered}
v_1=(6\;7)(12\;13),\quad v_2=1,\quad
v_3=(2\;3)(4\;6\;5)(7\;8\;9)(10\;12\;11)(13\;14),\\
v_4=(2\;3\;6\;5\;4)(7\;8\;9\;12\;11\;10)(13\;14),
\end{gathered}
\end{equation}
which have length $2$, $0$, $8$ and $10$, respectively.  The cubicles
are labelled accordingly in Figure \ref{figure;peacock}; the cubicle
$\ti\t_1=\t_+\cap\ti\t$ is shown in dark grey.  We determined $W\rel$
by picking elements in the interior of each cubicle, namely
$\ti\xi_1=3\ti\alpha\spcheck_1+2\ti\alpha\spcheck_2$,
$\ti\xi_2=6\ti\alpha\spcheck_1+5\ti\alpha\spcheck_2$,
$\ti\xi_3=5\ti\alpha\spcheck_1+6\ti\alpha\spcheck_2$, and
$\ti\xi_4=2\ti\alpha\spcheck_1+3\ti\alpha\spcheck_2$, and reflecting
them into the positive chamber $\t_+$.

This example disabused us of some overoptimistic notions we
entertained.  For instance, the subgroup $f\bigl(\SU(3)\bigr)$
consists of highly singular points in $\U(15)$.  Moreover, the small
Cartan $\ti\t$ has quite high codimension in $\bar\lie c$.  Notice
also that the cones $\ti\t_v$ come in two different shapes.
Furthermore, $\ti\alpha_2$ is not in the image of $S$, so $\ti S$ is
not a subset of $f^*(S)$, although $\ti R_+$ is a subset of $f^*(R_+)$
by Lemma \ref{lemma;face}.  Lastly, the relative Weyl elements do not
preserve the root system $\bar R$, so $W\rel$ is not a subset of $\ca
N_W(\bar C)$.
\end{example}

\section{Main results}\label{section;schubert}

As before, let $\ti K$ and $K$ be compact connected Lie groups and
$f\colon\ti K\to K$ a homomorphism with finite kernel.  In Section
\ref{subsection;inequal} we state our main result, a complete set of
inequalities for the moment cone of the cotangent bundle $T^*K$,
considered as a $\ti K\times K$-manifold.  We deduce from this
complete sets of inequalities for the moment cone of $T^*(\ti
K\backslash K)$ (considered as a $K$-manifold) and for the moment
polytope of every coadjoint orbit of $K$ (considered as a $\ti
K$-manifold).  As explained in the introduction, this is equivalent to
an asymptotic result on the decomposition of an irreducible $K$-module
into irreducible $\ti K$-modules.  This and other corollaries are
stated in Section \ref{subsection;corollaries}.  The inequalities we
obtain are in general highly overdetermined.  In Section
\ref{subsection;scalar} we explain how to prune them to a more
manageable set of inequalities.  We defer the proofs to Section
\ref{section;semistable}.  Finally, in Section
\ref{subsection;duality} we discuss the ``self-duality'' of the moment
cone of $T^*K$ and explain it in terms of the action of the longest
elements of $W$ and $\tW$.

\subsection{Schubert cycles and the moment cone}
\label{subsection;inequal}

Recall that $T$ and $\ti T$ denote maximal tori of $K$, resp.\ $\ti
K$, such that $\ti T=T\cap\ti K$.  Let $R\subset\t^*$ and $\ti
R\subset\ti\t^*$ be the respective root systems and $W$ and $\tW$
the associated Weyl groups.  Choose bases (sets of simple roots) $S$
in $R$ and $\ti S$ in $\ti R$ such that the corresponding (closed)
positive Weyl chambers $\t_+\subset\t$ and $\ti\t_+\subset\ti\t$ are
compatible in the sense of Definition \ref{definition;compatible}.
Denote the corresponding dual Weyl chambers in $\t^*$ and $\ti\t^*$ by
$\t^*_+$, resp.\ $\ti\t^*_+$.  Consider the symplectic manifold $T^*K$
on which $\ti K$ acts by left multiplication and $K$ by right
multiplication.  Let us identify $T^*K$ with $K\times\k^*$ by means of
left translations; then the action of $\ti K\times K$ is given by
$$
(\ti g,g)\cdot(k,\xi) =\bigl(f(\ti g)kg\inv,g\xi\bigr).
$$
With respect to the standard symplectic structure on $T^*K$ this
action has a moment map $\Phi =\Phi_{\ti K}\times\Phi_K\colon
T^*K\to\ti\k^*\times\k^*$ given by
\begin{equation}\label{equation;moment}
\Phi_{\ti K}(k,\xi) =f^*(k\xi),\qquad\Phi_K(k,\xi) =-\xi.
\end{equation}
According to \cite[Theorem 7.6]{sjamaar;convexity} the set
$$
\Delta(T^*K)=\Phi(T^*K)\cap(\ti\t^*_+\times\t^*_+)
$$
is a rational convex polyhedral cone, called the \emph{moment cone} of
$T^*K$.  Our main result is a set of inequalities describing this
cone.

\begin{theorem}\label{theorem;cone}
Let $(\tlambda,\lambda)\in\ti\t^*_+\times\t^*_+$.  Then
$(\tlambda,\lambda)\in\Delta(T^*K)$ if and only if
\begin{equation}\label{equation;inequal}
\ti w\inv\tlambda\in f^*(-w\inv\lambda+v\ca C)
\end{equation}
for all triples $(\ti w,w,v)\in\tW\times W\times W\rel$ such that
\begin{equation}\label{equation;nonzero}
\phi^*(v\sigma_{wv})(\ti c_{\ti w_0\ti w})\ne0.
\end{equation}
\end{theorem}

The proof is in Section \ref{section;semistable}.  First we explain
the statement and point out a number of consequences.  Recall that
$\ca C\subset\t^*$ denotes the root cone, i.e.\ the cone spanned by
the positive roots $R_+$.  The relative Weyl set $W\rel$ is the subset
(not subgroup) of the Weyl group $W$ defined in
\ref{definition;relative}.  For each $w\in W$ and $v\in W\rel$ the set
$f^*(-w\inv\lambda+v\ca C)$ is a polyhedral cone in $\ti\t^*$ with
apex $f^*(-w\inv\lambda)$, and so \eqref{equation;inequal} represents
a finite number of linear inequalities.  For instance, if $\ti\t_+$ is
contained in $\t_+$, then by Proposition \ref{proposition;face},
$W\rel=\{1\}$ and $f^*(\ca C)=\ti\ca C$, the root cone of $\ti K$.

The condition \eqref{equation;nonzero} is of a homological nature.
Let $X=K/T$ be the flag variety of $K$.  The Weyl group $W$ acts on
the homogeneous space $X$ by right multiplication, $w(kT)=kw\inv T$.
(We denote a Weyl group element and any of its representatives in $\ca
N_K(T)$ by the same letter as long as the formulas do not depend on
the choice of the representative.)  The induced action on a homology
or cohomology class $c$ is written as $wc$.

Recall that the homology of $X$ is torsion-free and has a basis given
by the Schubert classes $c_w\in H_{2l(w)}(X,\Z)$, where $w$ ranges
over $W$ and $l(w)$ is the length of $w$.  (The definition is reviewed
in Appendix \ref{section;flag}.)  If $w_0$ is the longest Weyl group
element, then $c_{w_0}$ is equal to the fundamental class $[X]$.  The
homology being torsion-free, we can identify the cohomology with the
dual abelian group $\Hom_\Z\bigl(H_\bu(X,\Z),\Z\bigr)$.  The
cohomology class $\sigma_w\in H^{2l(w)}(X,\Z)$ is defined by
$\sigma_w(c) =c\cdot c_{w_0w}$ for $c$ in $H_{2l(w)}(X,\Z)$ and
$\sigma_w(c)=0$ if $c$ has degree different from $2l(w)$.  Recall that
if $l(w)+l(w')\le l(w_0)$, then
\begin{equation}\label{equation;poincare}
c_w\cdot c_{w'} =\delta_{w,w_0w'}\,c_1,
\end{equation}
where the dot denotes the intersection product.  This implies that the
basis $\{\,\sigma_w\mid w\in W\,\}$ of $H^\bu(X,\Z)$ is dual to the
basis $\{\,c_w\mid w\in W\,\}$ of $H_\bu(X,\Z)$.

Finally, the map $\phi$ is the embedding of the flag variety $\ti
X=\ti K/\ti T$ into the flag variety $X=K/T$ which is induced by the
map $f\colon\ti K\to K$.  It is injective even if $f$ is not.

The Schubert bases depend on the choice of the dominant chambers.
However, it is not hard to work out the formula for a change of basis
(see e.g.\ Remark \ref{remark;basis}) and to show that a different
choice of compatible chambers $\ti\t_+$ and $\t_+$ leads up to a
relabelling to the same set of inequalities.

Theorem \ref{theorem;cone} makes the inequalities of the moment cone
computable in practice once we can

\begin{enumerate}
\item\label{item;relative}
determine the relative Weyl set $W\rel$,
\item\label{item;span}
write down inequalities for the codimension-one faces of the cones
$f^*(v\ca C)$, in other words, determine the \emph{rays} that span
their dual cones,
\item\label{item;action}
write the action of the elements of $W\rel$ on the homology in terms
of the Schubert basis,
\item\label{item;entry}
find the nonzero entries of the matrix of $\phi^*$ with respect to the
Schubert bases.
\end{enumerate}

Unfortunately, we do not know of a general solution to these
combinatorial problems except the third.  (Formulas for the Weyl group
action on the homology in terms of the Schubert basis were given by
Bernstein et al.\ \cite{bernstein-gelfand-gelfand;schubert} and
Demazure \cite{demazure;desingularisation}.  These are reviewed in
Appendix \ref{section;flag}.)  In many interesting examples, however,
they can be solved explicitly.  Some properties that can be helpful in
computing $W\rel$ and the cones $f^*(v\ca C)$ are discussed in Section
\ref{subsection;relative}.  Problem \eqref{item;entry} is discussed in
Appendix \ref{section;flag}.

\subsection{Alternative formulations and corollaries}
\label{subsection;corollaries}

Let $\lambda\in\t^*_+$ and consider the coadjoint orbit $\ca
O_\lambda=K\lambda$ through $\lambda$.  It follows from
\eqref{equation;moment} that the symplectic quotient $\Phi_K\inv(\ca
O_\lambda)/K$ of $T^*K$ with respect to the right $K$-action is
isomorphic, as a Hamiltonian $\ti K$-manifold, to the \emph{dual}
orbit $\ca O_{-w_0\lambda}$.  Consequently the moment polytope
$\Delta(\ca O_{-w_0\lambda})$ of $\ca O_{-w_0\lambda}$ with respect to
the $\ti K$-action is equal to the horizontal slice
$$
\Delta(\ca O_{-w_0\lambda})
=\Delta(T^*K)\cap\bigl(\ti\t^*\times\{\lambda\}\bigr)
$$
of the moment cone $\Delta(T^*K)$.  Furthermore, it is evident that
$\tlambda\in\Delta(\ca O_\lambda)$ if and only if $-\ti
w_0\tlambda\in\Delta(\ca O_{-w_0\lambda})$.  Thus, after substituting
$\ti w\to\ti w_0w$, we see that Theorem \ref{theorem;cone} is
equivalent to the following statement.

\begin{theorem}\label{theorem;polytope}
Let $(\tlambda,\lambda)\in\ti\t^*_+\times\t^*_+$.  Then
$\tlambda\in\Delta(\ca O_\lambda)$ if and only if
\begin{equation}\label{equation;polyinequal}
\ti w\inv\tlambda\in f^*(w\inv\lambda-v\ca C)
\end{equation}
for all triples $(\ti w,w,v)\in\tW\times W\times W\rel$ such that
$\phi^*(v\sigma_{wv})(\ti c_{\ti w})\ne0$.
\qed
\end{theorem}

Because of the equivalence of problems \eqref{item;geometric} and
\eqref{item;asymptotic} discussed in the introduction, Theorem
\ref{theorem;polytope} is tantamount to an asymptotic statement about
irreducible representations.  Let $\Lambda=\ker({\exp}|_\t)$ be the
integral lattice in $\t$, $\Lambda^*=\Hom_\Z(\Lambda,\Z)$ the weight
lattice and $\Lambda^*_+=\Lambda^*\cap\t^*_+$ the monoid of dominant
weights.  For every dominant weight $\lambda$ let $V_\lambda$ denote
the irreducible $K$-module with highest weight $\lambda$.

\begin{theorem}\label{theorem;representation}
Let $(\tlambda,\lambda)\in\ti\Lambda^*_+\times\Lambda^*_+$.  Then
there exists a positive integer $n$ such that $\ti V_{n\tlambda}$
occurs in $V_{n\lambda}$ if and only if \eqref{equation;polyinequal}
holds for all $(\ti w,w,v)\in\tW\times W\times W\rel$ such that
$\phi^*(v\sigma_{wv})(\ti c_{\ti w})\ne0$.
\qed
\end{theorem}

Similarly, it follows from \eqref{equation;moment} that the symplectic
quotient $\ti K\backslash\Phi_{\ti K}\inv(0)$ of $T^*K$ with respect
to the left $\ti K$-action is symplectomorphic to the cotangent bundle
of the homogeneous space $\ti K\backslash K$.  Hence the moment cone
of $T^*(\ti K\backslash K)$ with respect to the residual $K$-action is
equal to the ``vertical slice''
$$
\Delta\bigl(T^*(\ti K\backslash K)\bigr) =\Delta(T^*
K)\cap\bigl(\{0\}\times\t^*\bigr)
$$
of the moment cone $\Delta(T^*K)$.  The following result is now
immediate from Theorem \ref{theorem;cone}.

\begin{theorem}\label{theorem;smallcone}
Let $\lambda\in\t^*_+$.  Then $\lambda\in\Delta\bigl(T^*(\ti
K\backslash K)\bigr)$ if and only if
\begin{equation}\label{equation;null}
0\in f^*(w\inv\lambda-v\ca C)
\end{equation}
for all pairs $(w,v)\in W\times W\rel$ such that
$\phi^*(v\sigma_{wv})\ne0$.
\qed
\end{theorem}

Like Theorem \ref{theorem;polytope}, this has a
representation-theoretic counterpart.

\begin{theorem}\label{theorem;invariant}
Let $\lambda\in\Lambda^*_+$.  Then there exists $n>0$ such that
$V_{n\lambda}$ contains a $\ti K$-invariant vector if and only if
\eqref{equation;null} holds for all $(w,v)\in W\times W\rel$ such that
$\phi^*(v\sigma_{wv})\ne0$.
\qed
\end{theorem}

Although we have presented Theorems \ref{theorem;smallcone} and
\ref{theorem;invariant} as corollaries of Theorem \ref{theorem;cone},
they are actually equivalent to it and play a part in the proof; see
Section \ref{subsection;proof}.

Suppose now that $\ti V_\tlambda$ occurs in $V_\lambda$.  Then the
weights in $\ti V_\tlambda$ are a subset of the weights of $V_\lambda$
(with respect to the torus $\ti T$).  The latter are of the form
$f^*(\lambda-\mu)$, where $\mu$ is a combination of positive roots of
$K$ with nonnegative integral coefficients.  In particular, $\tlambda$
itself is of this form.  From Theorem \ref{theorem;representation} we
thus obtain the following necessary condition for $V_\lambda$ to
contain $\ti V_\tlambda$.

\begin{theorem}\label{theorem;lattice}
Let $(\tlambda,\lambda)\in\ti\Lambda^*_+\times\Lambda^*_+$.  If $\ti
V_\tlambda$ occurs in $V_\lambda$\upn, then
$$
\ti w\inv\tlambda\in f^*(w\inv\lambda-v\ca C_\Z)
$$
for all triples $(\ti w,w,v)\in\tW\times W\times W\rel$ such that
$\phi^*(v\sigma_{wv})(\ti c_{\ti w})\ne0$.
\qed
\end{theorem}

Here we have written $\ca C_\Z$ for the monoid generated by $R_+$.  An
interesting question is for what pairs $K$, $\ti K$ this necessary
condition is sufficient.  

\begin{example}\label{example;saturate}
The answer to this question is clearly affirmative when $\ti K=T$, the
maximal torus of $K$.  The same is true for the diagonal embedding of
$\U(n)$ into $\U(n)\times\U(n)$, as was recently shown by Knutson and
Tao \cite{knutson-tao;honeycomb-model-I}.  However, it was pointed out
to us by Knop that the answer is negative for the diagonal embedding
$\G_2\to\G_2\times\G_2$: the module $V_{03}$ is contained in
$V_{30}\otimes V_{03}$, but $V_{01}$ is not contained in
$V_{10}\otimes V_{01}$.  Here $V_{kl}$ is the irreducible module with
highest weight $k\pi_1+l\pi_2$, $\pi_1$ and $\pi_2$ being the
fundamental weights of $\G_2$.  This can be checked using a result of
Littelmann \cite[\S3.8]{littelmann;generalization}.  See Brion
\cite{brion;general-faces} and Montagar \cite{montagard;faces} for
further references and related results.
\end{example}

We discuss two more reformulations of Theorem \ref{theorem;cone}.
First, for all $(\ti w,w)\in\tW\times W$ let $\ca C_{\ti w,w}$ be the
cone $\bigcap_vf^*(v\ca C)$, where the intersection is taken over all
$v\in W\rel$ such that \eqref{equation;nonzero} holds (and put $\ca
C_{\ti w,w} =\ti\t^*$ if there is no such $v$).  Theorem
\ref{theorem;cone} is equivalent to:
$(\tlambda,\lambda)\in\Delta(T^*K)$ if and only if
$(\tlambda,\lambda)\in\ti\t^*_+\times\t^*_+$ and
$$
\ti w\inv\tlambda+f^*(w\inv\lambda)\in\ca C_{\ti w,w}
$$
for all $(\ti w,w)$.

Next, let $\ca P\colon H^i(X,\Z)\to H_{2l(w_0)-i}(X,\Z)$ denote the
Poincar\'e duality map, $\ca P(\sigma)=\sigma\cap[X]$.  Then $\ca
P(\sigma_w) =c_{w_0w}$ by \eqref{equation;poincare}.  Let
$$
\phi^! =\ti\ca P\circ\phi^*\circ\ca P\inv\colon H_\bu(X,\Z)\to
H_{\bu+2(l(\ti w_0)-l(w_0))}(\ti X,\Z)
$$
be the wrong-way or Gysin homomorphism induced by $\phi$.  If $l(\ti
w)+l(wv)=l(\ti w_0)$, then $\deg\phi^*(v\sigma_{wv}) =\deg\ti c_{\ti
w_0\ti w}$, so
\begin{equation}\label{equation;classes}
\begin{split}
\phi^*(v\sigma_{wv})(\ti c_{\ti w_0\ti w})&
=\ti\ca P\bigl(\phi^*(v\sigma_{wv})\bigr)\cdot\ti c_{\ti w_0\ti w}\\
& =\phi^!\bigl(\ca P(v\sigma_{wv})\bigr)\cdot\ti c_{\ti w_0\ti w}\\
& =(-1)^{l(v)}\phi^!(v\ca P\sigma_{w v})\cdot\ti c_{\ti w_0\ti w}\\
& =(-1)^{l(v)}\phi^!(vc_{w_0wv})\cdot\ti c_{\ti w_0\ti w}.
\end{split}
\end{equation}
Here we have used that the mapping degree of the diffeomorphism of $X$
induced by $v$ is equal to $(-1)^{l(v)}$, which follows from Theorem
\ref{theorem;dbgg}.  Let us say that a Schubert class $c_w$ is
\emph{contained} in a homology class $c$ if it occurs with nonzero
coefficient in the expression for $c$ as a linear combination of
Schubert classes.  This amounts to $\sigma_w(c)\ne0$, or equivalently,
if $c$ is homogeneous, $c\cdot c_{w_0w}\ne0$ and $\deg c =2l(w)$.
Then according to \eqref{equation;classes}, condition
\eqref{equation;nonzero} is equivalent to
\begin{equation}\label{equation;contain}
\ti c_{\ti w}\quad\text{is contained in}\quad\phi^!(vc_{w_0 wv}).
\end{equation}
This means that there exist cycles on $\ti X$ and $X$ which are
homologous to $\ti c_{\ti w}$ and $vc_{w_0wv}$, respectively, and
which have a nontrivial transverse intersection.  Applying $\ti\ca
P\inv$ to both sides yields a cohomological version of
\eqref{equation;contain}:
$$
\ti\sigma_{\ti w_0\ti w}\quad\text{is contained
in}\quad\phi^*(v\sigma_{wv}).
$$
Observe also that the condition $\phi^*(v\sigma_{wv})(\ti c_{\ti
w})\ne0$ of Theorem \ref{theorem;polytope} is equivalent to
$$
\ti\sigma_{\ti w}\quad\text{is contained in}\quad\phi^*(v\sigma_{wv}).
$$
This reformulation of Theorem \ref{theorem;polytope} does not refer to
Poincar\'e duality or the longest Weyl group elements and therefore
makes sense, at least syntactically, for coadjoint orbits of certain
infinite-dimensional groups.

\subsection{Scalar inequalities}\label{subsection;scalar}

Many of the vector inequalities \eqref{equation;inequal} are
redundant.  In fact, each of them can be replaced by a single scalar
inequality, and in addition \eqref{equation;nonzero} can be reduced to
a much smaller set of homological conditions.  To do so we need to
recall some more results from
\cite{bernstein-gelfand-gelfand;schubert}.

Let $\s$ be a (closed) face of the Weyl chamber $\t_+$.  The
centralizer of $\s$ is equal to the centralizer of any point in its
(relative) interior $\s^\circ$ and is denoted by $\ca Z_K(\s)$ or
$K_\s$.  It is a connected subgroup of $K$ and contains the maximal
torus $T$.  Its Weyl group $\ca N_{K_\s}(T)/T$ is denoted by $W_\s$,
and the partial flag variety $K/K_\s$ by $X_\s$.  The homology of
$X_\s$ has a free basis consisting of Schubert classes $c^\s_w\in
H_{2l(w)}(X_\s,\Z)$, where $w$ ranges over the set $W^\s$ defined by
\begin{equation}\label{equation;minimal}
W^\s =\bigl\{\,w\in W\bigm|l(ww')\ge l(w)\quad\text{for all $w'\in
W_\s$}\,\bigr\}.
\end{equation}
Every coset $wW_\s$ contains a unique element in $W^\s$, called its
\emph{shortest representative}.  The canonical projection $X\to X_\s$
induces a surjection $H_\bu(X,\Z)\to H_\bu(X_\s,\Z)$, which maps $c_w$
to $c_w^\s$ if $w\in W^\s$ and to $0$ if $w\not\in W^\s$.  It also
induces an embedding of $H^\bu(X_\s,\Z)$ onto the subspace of
$W_\s$-invariants in $H^\bu(X,\Z)$.  A basis of $H^\bu(X,\Z)^{W_\s}$
is given by the Schubert cocycles $\{\,\sigma_w\mid w\in W^\s\,\}$.
Via the identification $H^\bu(X_\s,\Z)\cong H^\bu(X,\Z)^{W_\s}$ this
basis is dual to $\{\,c^\s_w\mid w\in W^\s\,\}$.  The Poincar\'e dual
of $\sigma_w$ is $c^\s_u$, where $u$ is the shortest representative of
$w_0wW_\s$.

Consider any point $\ti\chi$ in $\ti\t_+$ and let $\ti\s$ be the face
in $\ti\t$ which contains $\ti\chi$ in its interior.  Recall that
$\ti\t_v$ denotes the cone $v\t_+\cap\ti\t$.  Choose $v\in W\rel$ such
that $\ti\chi\in\ti\t_v$, in other words $v\inv\ti\chi\in\t_+$, and
let $\s$ be the face of $\t_+$ which contains $v\inv\ti\chi$ in its
interior.  Then clearly $\ti K_{\ti\s}=vK_\s v\inv\cap\ti K$, so we
can define an equivariant embedding
\begin{equation}\label{equation;embed}
\phi_{\ti\s,v}\colon\ti X_{\ti\s}\to X_\s
\end{equation}
by sending $\ti\pi_{\ti\s}(1)$ to $\pi_\s(v)$.  (Here $\pi_\s$ denotes
the quotient map $K\to X_\s$.)  Now select nonzero rational elements
$\ti\chi_1$, $\ti\chi_2,\dots$, $\ti\chi_n\in\ti\t$ such that each of
the subcones $\ti\t_v$ of $\ti\t_+$ is spanned by an appropriate
subcollection.  If $\ti\t\cap\lie c=0$, this can be accomplished by
selecting one nonzero element on every ray (one-dimensional face)
occurring in the cubicle subdivision \eqref{equation;subdivision}.
For $k=1$, $2,\dots$, $n$, pick $v_k\in W\rel$ such that
$\ti\chi_k\in\ti\t_{v_k}$.  Let $\ti\s_k$ and $\s_k$ be the faces of
$\ti\t$ and $\t$ such that $\ti\chi_k\in\ti\s^\circ_k$ and
$v_k\inv\ti\chi_k\in\s^\circ_k$, respectively.  For brevity put
$W_{\s_k}=W_k$, $W^{\s_k}=W^k$, and $\phi_{\ti\s_k,v_k}=\phi_k$.

\begin{theorem}\label{theorem;scalar}
Let $(\tlambda,\lambda)\in\ti\t^*_+\times\t^*_+$.  Then
$(\tlambda,\lambda)\in\Delta(T^*K)$ if and only if
$$
\bigl\langle\ti w\inv\tlambda
+f^*(w\inv\lambda),\ti\chi_k\bigr\rangle\ge0
$$
for $k=1$\upn, $2,\dots$\upn, $n$ and for all $(\ti w,w)\in\tW^k\times
W^k$ such that $\ti\sigma_{\ti u_k}$ is contained in
$\phi_k^*(\sigma_{u_k})$.  Here $\ti u_k\in\tW^k$ and $u_k\in W^k$
denote the shortest representatives of $\ti w_0\ti w\tW_k$ and
$wv_kW_k$\upn, respectively.
\end{theorem}

For instance, if $\ti\t_+$ is contained in $\t_+$ and $\ti K$ is
semisimple, then we can choose the $\chi_k$ to be the fundamental
coweights.  Then the cohomological conditions in this theorem involve
only the minimal partial flag varieties (``Grassmannians'') of $\ti K$
and certain partial flag varieties of $K$, and hence the number of
inequalities is considerably smaller than in Theorem
\ref{theorem;cone}.  In Section \ref{section;examples} we shall see,
however, that even the system of inequalities of Theorem
\ref{theorem;scalar} is often overdetermined.

\subsection{Duality}\label{subsection;duality}

The momentum map on $T^*K$ has the property that
$\Phi(k,-\xi)=-\Phi(k,\xi)$.  This implies immediately that
$\Delta(T^*K)$ is stable under the \emph{duality involution}, which
sends $(\tlambda,\lambda)$ to $(-\ti w_0\tlambda,-w_0\lambda)$.
Consider the inequality \eqref{equation;inequal} associated with a
triple $(\ti w,w,v)\in\tW\times W\times W\rel$.  Substituting
$\tlambda\to-\ti w_0\tlambda$ and $\lambda\to-w_0\lambda$ and
multiplying both sides by $\ti w_0$ we obtain
\begin{equation}\label{equation;winequal}
\ti w_0\ti w\inv\ti w_0\tlambda\in-\ti w_0f^*(w\inv w_0\lambda+v\ca
C).
\end{equation}

This motivates the following definition.  Let
$j\colon\tW\hookrightarrow W$ be the lifting homomorphism associated
with the chambers $\t_+$ and $\ti\t_+$.  (Cf.\ Remark
\ref{remark;choice}.)

\begin{definition}
The \emph{dual} of $w\in W$ is $w^*=j(\ti w_0)ww_0$.  The \emph{dual}
of $\ti w\in\tW$ is $\ti w^*=\ti w_0w\ti w_0$.
\end{definition}

Observe that $(w^*)^*=w$ and $(\ti w^*)^*=\ti w$.  However, the
duality map on $W$ is in general not a homomorphism, nor does it
preserve the length function on $W$ or map the subgroup $\tW$ into
itself.

\begin{lemma}\label{lemma;dual}
The set of compatible Weyl group elements $W\com\subset W$ is stable
under the duality map.  Furthermore $\ti\t_{w^*}=-\ti w_0\ti\t_w$ and
$f^*(v^*\ca C)=-\ti w_0f^*(v\ca C)$ for $w\in W\com$.  Hence the
subdivision into cubicles \eqref{equation;subdivision} is stable under
the duality involution of $\ti\t_+$.
\end{lemma}

\begin{proof}
Using the definition of $W\com$ and the fact that $w_0\t_+=-\t_+$ we
find
\begin{align*}
w\in W\com&\iff w\t_+\cap\ti\t_+^\circ\ne\emptyset\\
%
&\iff j(\ti w_0)ww_0\t_+\cap\ti\t_+^\circ\ne\emptyset\\
&\iff w^*\in W\com.
\end{align*}
By the same token, $\ti\xi\in\ti\t_{w^*}$ if and only if $\ti\xi
=j(\ti w_0)ww_0\xi$ for some $\xi\in\t_+$, which is equivalent to
$-\ti w_0\ti\xi =w(-w_0\xi)\in w\t_+\cap\ti\t=\ti\t_w$.  In other
words, $\ti\t_{w^*}=-\ti w_0\ti\t_w$.  The equality $f^*(v^*\ca
C)=-\ti w_0f^*(v\ca C)$ now follows from Lemma \ref{lemma;dualcone},
and the last statement is obvious.
\end{proof}

\begin{example}
Consider the map $f\colon\SU(3)\to\U(15)$ studied in Example
\ref{example;su3}.  By \eqref{equation;j},
\begin{gather*}
j(\ti w_0) =j(\ti s_1\ti s_2\ti s_1)
=(1\;15)(2\;13)(3\;14)(4\;10)(5\;11)(6\;12),\\
\intertext{and therefore}
v_1^*=(2\;3\;6\;4)(7\;9\;12\;10)(13\;14),\\
v_2^*=1^*=(2\;3)(4\;6)(7\;9)(10\;12)(13\;14).
\end{gather*}
Comparing with \eqref{equation;relative} we see that the duality map
does not map $W\rel$ into itself, but for each $v\in W\rel$,
$v^*\in\barW v'$ for some $v'\in W\rel$.
\end{example}

Together with \eqref{equation;winequal} Lemma \ref{lemma;dual} yields
\begin{equation}\label{equation;dualinequal}
(\ti w\inv)^*\tlambda\in f^*\bigl(-(w\inv)^*\lambda+v^*\ca C\bigr),
\end{equation}
which is called the \emph{dual} inequality associated with
\eqref{equation;inequal}.  It can be obtained from
\eqref{equation;inequal} by replacing the triple $(\ti w,w,v)$ with
$\bigl(\ti w^*,((w\inv)^*)\inv,v^*\bigr)$.  We assert that the dual
inequality arises from a homological condition analogous to
\eqref{equation;nonzero}, namely
\begin{equation}\label{equation;oppositenonzero}
\phi^*(v^*\sigma_{((w\inv)^*)\inv v^*})(\ti c_{\ti w_0\ti w^*})\ne0.
\end{equation}
Indeed,
\begin{align*}
\phi^*(v^*\sigma_{((w\inv)^*)\inv v^*})(\ti c_{\ti w_0\ti w^*})
&=\phi^*\bigl(j(\ti w_0)vw_0\sigma_{w_0wvw_0}\bigr)(\ti c_{\ti w\ti
w_0})\\
&=\phi^*(vw_0\sigma_{w_0wvw_0})(\ti w_0\ti c_{\ti w\ti w_0})\\
&=(-1)^{l(wv)+l(\ti w_0)-l(\ti w)}\phi^*(v\sigma_{wv})(\ti c_{\ti
w_0\ti w}),
\end{align*}
where we have used the $\tW$-equivariance of $\phi$ (see Corollary
\ref{corollary;lift}) and Lemma \ref{lemma;longact} below.  This shows
that \eqref{equation;nonzero} implies
\eqref{equation;oppositenonzero}, which explains the appearance of the
dual inequality \eqref{equation;dualinequal}.

\begin{lemma}\label{lemma;longact}
$w_0c_w=(-1)^{l(w)}c_{w_0ww_0}$ and
$w_0\sigma_w=(-1)^{l(w)}\sigma_{w_0ww_0}$.
\end{lemma}

\begin{proof}
The second equality is immediate from the first.  By Remark
\ref{remark;basis}, to prove the first equality it is enough to show
that $c^{w_0}_w=(-1)^{l(w)}c_w$.  We use the notation of that remark.
The \emph{opposite} or \emph{dual} Borel subgroup is
$B^{w_0}=w_0Bw_0$.  Let $\theta\colon G\to G$ be complex conjugation
with respect to the real form $K$; then $\theta$ maps $B$ to its
opposite and therefore induces an antiholomorphic map $G/B\to
G/B^{w_0}$, which we also denote by $\theta$.  If $Y$ is a complex
submanifold of $G/B$, then so is $\theta(Y)$, but $\theta$ changes the
orientation if the complex dimension of $Y$ is odd.  It is clear that
the diagram
$$
\begin{CD}
X@>\tau>>G/B\\
@V{\id_X}VV@VV{\theta}V\\
X@>{\tau_{w_0}}>>G/B^{w_0}
\end{CD}
$$
commutes.  Moreover, $\theta(X_w^\circ)=(X^{w_0}_w)^\circ$ and
therefore
$$
(\tau_{w_0})_*(c_w) =\theta_*\tau_*(c_w) =\theta_*\bigl([X_w]\bigr)
=(-1)^{l(w)}[X^{w_0}_w] =(-1)^{l(w)}(\tau_{w_0})_*c^{w_0}_w,
$$
which shows that $c^{w_0}_w=(-1)^{l(w)}c_w$.
\end{proof} 

\section{Semistability}\label{section;semistable}

This section contains the proofs of Theorems \ref{theorem;cone} and
\ref{theorem;scalar}.  In Sections \ref{subsection;measure} and
\ref{subsection;flag} we review the Hilbert-Mumford criterion and
calculate Mumford's numerical measure of instability for flag
varieties.  We then finish the proof in Section
\ref{subsection;proof}.  The underlying idea is that an integral
coadjoint orbit is a complex projective variety in a natural way, so
we can detect points in its momentum polytope by determining the
semistable set with respect to various ample line bundles.  This is
achieved by means of the Hilbert-Mumford criterion, and the upshot is
that the semistable set is nonempty if and only if a certain translate
of a Schubert cell in the flag variety $X$ intersects the small flag
variety $\ti X$.  By a transversality argument this leads to
inequalities in terms of intersections of Schubert cycles.

Throughout this section the complexifications of the groups $K$ and
$\ti K$ are denoted by $G$ and $\ti G$, respectively.

\subsection{Measure of instability}\label{subsection;measure}

The following material is extracted from Mumford's book
\cite{mumford-fogarty-kirwan;geometric}.  Let $Y$ be a complex
projective $G$-variety (where $G=K^\C$) and let $\ca L$ be an ample
$G$-equivariant line bundle over $Y$.  Recall that $y\in Y$ is called
\emph{semistable} with respect to $\ca L$ if for some $n>0$ there
exists an invariant section $s$ of $\ca L^n$ such that $s(y)\ne0$.  A
point is \emph{unstable} if it is not semistable.  The set of
semistable points is denoted by $Y\sst$.  It is clearly Zariski-open
and $G$-stable.

Let $y$ be any point in $Y$ and let $\chi\in\Hom(\C^\times,G)$ be an
algebraic one-parameter subgroup of $G$.  Mumford's \emph{numerical
measure of instability} is the integer $m^{\ca L}(y,\chi)$ determined
as follows.  Consider $y_0=\lim_{t\to0}\chi(t)y$, which exists because
$Y$ is projective.  It is a fixed point for the action of $\chi$ and
so $\chi$ acts on the fibre $\ca L_{y_0}$.  Let $r$ be the unique
integer such that $\chi(t)\,l=t^rl$ for $l\in\ca L_{y_0}$.  Then
$m^{\ca L}(y,\chi)=r$.  (This sign convention is opposite to
Mumford's, but agrees with our sign convention for the moment map.)

\begin{theorem}[Hilbert-Mumford criterion]\label{theorem;hm}
A point $y\in Y$ is semistable if and only if $m^{\ca L}(y,\chi)\le0$
for all one-parameter subgroups $\chi$.
\end{theorem}

We will frequently use the following elementary properties of $m^{\ca
L}(y,\chi)$.

\begin{proposition}\label{proposition;indispensable}
\begin{enumerate}
\item\label{item;equivariant}
$m^{\ca L}(gy,g\chi g\inv)=m^{\ca L}(y,\chi)$ for all $g\in G$.
\item\label{item;invariant}
$m^{\ca L}(py,\chi)=m^{\ca L}(y,\chi)$ for all $p\in P_\chi$.
\item\label{item;additive}
If $(Y_1,\ca L_1)$ and $(Y_2,\ca L_2)$ are two linearized
$G$-varieties\upn, then
$$
m^{\ca L_1\boxtimes\ca L_2}\bigl((y_1,y_2),\chi\bigr) =m^{\ca
L_1}(y_1,\chi)+m^{\ca L_2}(y_2,\chi)
$$
for $(y_1,y_2)\in Y_1\times Y_2$.
\end{enumerate}
\end{proposition}

Here $P_\chi$ denotes the parabolic subgroup associated with $\chi$,
which is defined in \eqref{equation;parabolic}.  We identify
$\Hom(\C^\times,H)$, the set of algebraic one-parameter subgroups of
the complex Cartan $H=T^\C$, with the lattice $\Lambda\subset\t$ by
assigning to each one-parameter subgroup $\chi$ its infinitesimal
generator $d\chi(1)$ times $2\pi i$.  We call $\chi$ \emph{dominant}
if it is contained in $\Lambda_+=\Lambda\cap\t_+$ and denote by
$\langle\xi,\chi\rangle$ the natural pairing between $\chi$ and any
$\xi\in\g^*$.  According to
\ref{proposition;indispensable}\eqref{item;equivariant}, to calculate
$m^{\ca L}(y,\chi)$ for arbitrary one-parameter subgroups $\chi$ of
$G$ it suffices to calculate $m^{\ca L}(gy,\chi)$ for arbitrary $g\in
G$ and dominant one-parameter subgroups $\chi$ of $H$.

\subsection{Instability on flag varieties}\label{subsection;flag}

The following discussion generalizes
\cite[\S4.4]{mumford-fogarty-kirwan;geometric}.  Take $Y$ to be the
flag variety $X=K/T$.  We identify $X$ with the complex homogeneous
space $G/B$ and consider the homogeneous line bundle
$$
\ca L_\lambda=G\times^{B}\C_\lambda,
$$
where $\lambda$ is a strictly dominant weight and $\C_\lambda$ is the
one-dimensional representation of $T$ defined by the weight $\lambda$
(extended to $B$ by letting $H=T^\C$ act holomorphically and by
letting the unipotent radical of $B$ act trivially).  Instead of
$m^{\ca L_\lambda}$ we write $m^\lambda$.  Consider an arbitrary $x\in
X$.  Let $\pi\colon G\to X$ be the quotient map and write $x=\pi(g)$.
We assert that if $h\in G$ and $\chi$ is dominant,
\begin{equation}\label{equation;flag}
m^\lambda(x,h\chi h\inv)=m^\lambda\bigl(h\inv\pi(g),\chi\bigr)
=\langle w\inv\lambda,\chi\rangle,
\end{equation}
where $w$ is the unique Weyl group element for which $\pi(h)$ is in
the translated cell $gX^\circ_w$.  Here we use the notation
$X^\circ_w$ for the Bruhat cell $BwB/B\subset X$ corresponding to
$w\in W$.  Indeed, $\pi(h)\in gX^\circ_w$ is equivalent to
$h\inv\pi(g)$ being in the $B$-orbit of $\pi(w\inv)$.  Since $\chi$ is
dominant, $P_\chi$ contains $B$, so
$$
m^\lambda\bigl(h\inv\pi(g),\chi\bigr)
=m^\lambda\bigl(\pi(w\inv),\chi\bigr)
$$
by \ref{proposition;indispensable}\eqref{item;invariant}.  Because
$\pi(w\inv)$ is fixed under $H$, we can calculate the right-hand side
by considering the action of $\chi$ on the fibre of $\ca
L_\lambda=G\times^B\C_\lambda$ over $\pi(w\inv)$:
$$
\chi(t)[w\inv,z]=[\chi(t)w\inv,z]
=\bigl[w\inv,\bigl(w\chi(t)w\inv\bigr)z\bigr] =\bigl[w\inv,t^{\langle
w\inv\lambda,\chi\rangle}z\bigr].
$$
This proves $m^\lambda\bigl(\pi(w\inv),\chi\bigr) =\langle
w\inv\lambda,\chi\rangle$ and hence \eqref{equation;flag}.

Evidently, there are no $G$-semistable points on $X$, but the
situation becomes interesting when we restrict the action of $G$ to
the subgroup $\ti G=\ti K^\C$.  Again we need only calculate
$m^\lambda(x,\ti g\ti\chi\ti g\inv)=m^\lambda\bigl(\ti
g\inv\pi(g),\ti\chi\bigr)$ for arbitrary $\ti g\in\ti G$ and dominant
$\ti\chi\in\Hom(\C^\times,\ti H)$.  Now $\ti\chi$ is not necessarily
dominant for $G$, but according to \eqref{equation;subdivision} we can
choose $v\in W\rel$ such that $v\inv\ti\chi$ is dominant for $G$.
Then $m^\lambda\bigl(\ti g\inv\pi(g),\ti\chi\bigr)
=m^\lambda\bigl(v\inv\ti g\inv\pi(g),v\inv\ti\chi\bigr)$, which can be
computed by \eqref{equation;flag}.  The upshot is:
\begin{equation}\label{equation;bigflag}
m^\lambda(x,\ti g\ti\chi\ti g\inv)= m^\lambda\bigl(\ti
g\inv\pi(g),\ti\chi\bigr) =\bigl\langle
f^*(w\inv\lambda),\ti\chi\bigr\rangle,
\end{equation}
where $w\in W$ is determined by the condition that $\ti gv\in gBwvB$,
i.e.\ $\ti g\pi(v)\in gX^\circ_{wv}$.  By Proposition
\ref{proposition;plusface}\eqref{item;borel} the stabilizer of
$\pi(v)\in X$ with respect to the $\ti G$-action is the Borel $\ti B$.
Let $\phi_v\colon\ti X\to X$ be the $\ti G$-equivariant embedding
which sends $\ti\pi(1)$ to $\pi(v)$.  Using the Hilbert-Mumford
criterion we obtain the following result.

\begin{proposition}\label{proposition;semi}
Let $\lambda$ be a strictly dominant weight.  A point $\pi(g)$ in
$X\cong G/B$ is semistable with respect to the bundle $\ca L_\lambda$
and the action of the subgroup $\ti G$ if and only if
\begin{equation}\label{equation;less}
\bigl\langle f^*(w\inv\lambda),\ti\chi\bigr\rangle\le0
\end{equation}
for all $(w, v)\in W\times W\rel$ and $\ti\chi\in\ti\Lambda$ such that
$\ti\chi\in\ti\Lambda_+\cap v\Lambda_+$ and $g
X^\circ_{wv}\cap\phi_v(\ti X)$ is nonempty.
\qed
\end{proposition}

If $\lambda$ is dominant but not strictly dominant, then it does not
define an ample bundle on $X$, but on the partial flag variety
$X_\s\cong G/P$, where $\s$ is the face of $\t_+$ which contains
$\lambda$ in its interior, and $P=P_\s$, the parabolic associated with
$\s$.  The analogue of \eqref{equation;bigflag} is: for all $g$ and
dominant $\ti\chi$
$$
m^\lambda\bigl(\ti g\inv\pi_\s(g),\ti\chi\bigr) =\bigl\langle
f^*(w\inv\lambda),\ti\chi\bigr\rangle,
$$
where $w$ is any element of $W$ (unique modulo \emph{left}
multiplication by $W_\s$) such that $\ti gv\in gPwvB$, i.e.\
$\phi_v\bigl(\ti\pi(\ti g)\bigr)\in g\bigl(P\pi(wv)\bigr)$.  (Note:
this condition involves a $P$-orbit in $G/B$, not a $B$-orbit in
$G/P$.)  By analogy with Proposition \ref{proposition;semi} we obtain
the following.

\begin{proposition}\label{proposition;partialsemi}
Let $\lambda\in\Lambda^*\cap\s^\circ$.  A point $\pi_\s(g)$ in
$X_\s\cong G/P$ is semistable with respect to the bundle $\ca
L_\lambda$ and the action of the subgroup $\ti G$ if and only if
\eqref{equation;less} holds for all $(w, v)\in W\times W\rel$ and
$\ti\chi\in\ti\Lambda$ such that $\ti\chi\in\ti\Lambda_+\cap
v\Lambda_+$ and $g\bigl(P\pi(wv)\bigr)\cap\phi_v(\ti X)$ is nonempty.
\qed
\end{proposition}

It is in fact not necessary to verify \eqref{equation;less} for all
dominant one-parameter subgroups of $\ti H$.  Choose a collection of
dominant one-parameter subgroups $\ti\chi_1$, $\ti\chi_2,\dots$,
$\ti\chi_n$ of $\ti H$ such that each of the cubicles $\ti\t_v$ in
$\ti\t_+$ is spanned by an appropriate subcollection.  Then every
dominant one-parameter subgroup $\ti\chi$ can be written as
\begin{equation}\label{equation;generator}
\ti\chi=\frac1{a}\sum_{k=1}^na_k\ti\chi_k,
\end{equation}
where $a$ and $a_k$ are integers, $a>0$, $a_k\ge0$, and $a_k=0$ if
$\ti\chi_k$ is not in the cubicle containing $\ti\chi$.  Now let
$\lambda\in\Lambda^*_+$ be strictly dominant.  It follows from
\eqref{equation;bigflag} and \eqref{equation;generator} that
$\pi(g)\in X$ is semistable for the action of $\ti G$ if and only if
$m^\lambda\bigl(\ti g\inv\pi(g),\ti\chi_k\bigr)\le0$ for all $k$ and
for all $\ti g\in\ti G$.  For $k=1$, $2,\dots$, $n$, fix $v_k\in
W\rel$ such that $\ti\chi_k\in\ti\t_{v_k}$ (i.e.\ $v_k\inv\ti\chi_k$
is dominant relative to $G$) and let $\ti\s_k$, resp.\ $\s_k$, be the
face in $\ti\t_+$, resp.\ $\t_+$, which is determined by
$\ti\chi_k\in\ti\s^\circ_k$, resp.\ $v\inv\ti\chi_k\in\s^\circ_k$.
For simplicity let us write $X_k=X_{\s_k}$, $P_k=P_{\s_k}$,
$W_{\s_k}=W_k$, $W^{\s_k}=W^k$, and $\ti\pi_{\ti\s_k}=\ti\pi_k$.  Let
$\phi_{\ti\s_k,v_k}=\phi_k$ be the embedding of $\ti X_k$ into $X_k$
defined in \eqref{equation;embed}.  Again using
\ref{proposition;indispensable}\eqref{item;invariant} we get a
stronger version of \eqref{equation;bigflag}:
$$
m^\lambda\bigl(\ti g\inv\pi(g),\ti\chi_k\bigr) =\bigl\langle
f^*(w\inv\lambda),\ti\chi_k\bigr\rangle,
$$
where $w\in W^k$ is determined by $\ti gv_k\in gBwv_kP_k$, in other
words $\phi_k\bigl(\ti\pi_k(\ti g)\bigr)\in g(X_k)^\circ_{wv_k}$.  The
conclusion is as follows.

\begin{proposition}\label{proposition;maximal}
Let $\lambda$ be a strictly dominant weight.  A point $\pi(g)$ in $X$
is semistable with respect to the bundle $\ca L_\lambda$ and the
action of the subgroup $\ti G$ if and only if
$$
\bigl\langle f^*(w\inv\lambda),\ti\chi_k\bigr\rangle\le0
$$
for $k=1$\upn, $2,\dots$\upn, $n$ and for all $w\in W^k$ such that
$g(X_k)^\circ_{u_k}\cap\phi_k(\ti X_k)$ is nonempty.  Here $u_k\in
W^k$ denotes the shortest representative of $wv_kW_k$.
\qed
\end{proposition}

We leave it to the reader to state a version for singular dominant
weights.

\subsection{Proof of Theorems \ref{theorem;cone} and
\ref{theorem;scalar}}\label{subsection;proof}

The proof of Theorem \ref{theorem;cone} is in two steps: first we show
that Theorem \ref{theorem;invariant} implies Theorem
\ref{theorem;representation}, and then we establish Theorem
\ref{theorem;representation}.  As we pointed out earlier, Theorem
\ref{theorem;representation} is equivalent to Theorem
\ref{theorem;cone}, so this will finish the proof.

\emph{Step}~1.  Suppose we knew Theorem \ref{theorem;invariant} was
true.  Let us apply it to the pair $(\ti K,\ti K\times K)$ and the
inclusion map $\id_{\ti K}\times f$.  Consider
$(\tlambda,\lambda)\in\ti\Lambda^*_+\times\Lambda^*_+$ and note that
$(\ti V_{\tlambda}\otimes V_{\lambda})^{\ti K}\ne\{0\}$ if and only if
$\ti V_{\tlambda}$ is a summand of $V_{\lambda}^*\cong
V_{-w_0\lambda}$.  Hence, by Theorem \ref{theorem;smallcone}, there
exists $n>0$ such that $\ti V_{n\tlambda}$ occurs in
$V_{-nw_0\lambda}$ if and only if
\begin{equation}\label{equation;otherinequal}
0\in(\id_{\ti K}\times f)^*\bigl((\ti
w\inv\tlambda,w\inv\lambda)-\ti\ca C\times v\ca C\bigr)
\end{equation}
for all $(\ti w,w,v)$ such that
\begin{equation}\label{equation;othernonzero}
(\id_{\ti X}\times\phi)^!(\ti c_{\ti w_0\ti w},vc_{w_0 wv})\ne0.
\end{equation}
Now $(\id_{\ti K}\times f)^*(\ti\ca C\times v\ca C) =\ti\ca C+f^*(v\ca
C)=f^*(v\ca C)$ by Proposition
\ref{proposition;plusface}\eqref{item;cone}, so
\eqref{equation;otherinequal} is equivalent to
\begin{equation}\label{equation;secondinequal}
\ti w\inv\tlambda\in f^*(-w\inv\lambda+v\ca C).
\end{equation}
Furthermore, \eqref{equation;othernonzero}
is equivalent to
\begin{equation}\label{equation;intersect}
\ti c_{\ti w_0\ti w}\cdot\phi^!(vc_{w_0wv})\ne0
\end{equation}
in $H_{2m}(\ti X,\Z)$, where $m=l(\ti w_0)-l(\ti w)-l(wv)$.  The
conclusion is that $\ti V_{n\tlambda}\subset V_{-nw_0\lambda}$ for
some $n>0$ if and only if \eqref{equation;secondinequal} holds for all
$\ti w$, $w$ and $v$ satisfying \eqref{equation;intersect}.

In \eqref{equation;intersect} it is actually sufficient to take into
account only those triples $(\ti w,w,v)$ for which the intersection is
zero-dimensional.  To see why, recall that the class $\ti c_{\ti
w_0\ti w}$ is represented by the Schubert variety $\ti X_{\ti w_0\ti
w}\subset\ti X$ and that $vc_{w_0wv}$ is represented by an algebraic
cycle $X_{w,v}$ which is a linear combination of Schubert varieties in
$X$.  Then \eqref{equation;intersect} means that for generic $\ti
g\in\ti K$ and $g\in K$ the translated cycles $\ti g\ti X_{\ti w_0\ti
w}$ and $gX_{w,v}$ intersect in an $m$-dimensional cycle on $\ti X$.
If $m>0$, this implies they intersect in a boundary component $\ti
X^\circ_{\ti w_0\ti u}$ of $\ti X_{\ti w_0\ti w}$ because of the fact
that $\ti X^\circ_{\ti w_0\ti w}$ is an affine variety.  In other
words, we have $\ti c_{\ti w_0\ti u}\cdot\phi^!(vc_{w_0wv})\ne0$ for
some $\ti u\suc\ti w$ (where $\suc$ denotes the Bruhat-Chevalley
order), and hence
\begin{equation}\label{equation;stronginequal}
\ti u\inv\tlambda\in f^*(-w\inv\lambda+v\ca C).
\end{equation}
Because $\tlambda$ is dominant, $\ti u\inv\tlambda\le\ti
w\inv\tlambda$, so if \eqref{equation;stronginequal} holds, then 
$$
\ti w\inv\tlambda\in\ti u\inv\tlambda+\ti\ca C\subset
f^*(-w\inv\lambda+v\ca C)+\ti\ca C= f^*(-w\inv\lambda+v\ca C).
$$
This shows that the inequality \eqref{equation;stronginequal} is
stronger than \eqref{equation;secondinequal}.  By induction we can
discard all inequalities \eqref{equation;intersect} except those for
which $m=0$, in which case they are equivalent to
$\phi^*(v\sigma_{wv})(\ti c_{\ti w_0\ti w})\ne0$.  Thus we see (after
replacing $\lambda$ by $-w_0\lambda$ and $w$ by $w_0w$) that Theorem
\ref{theorem;invariant} implies Theorem \ref{theorem;representation}.

\emph{Step}~2.  We now prove Theorem \ref{theorem;invariant}.
Consider $\lambda\in\Lambda^*_+$.  Let us assume first that $\lambda$
is strictly dominant.  Let $\ca L_\lambda$ the homogeneous ample line
bundle on $X\cong G/B$ introduced in Section \ref{subsection;flag}.
The Borel-Weil Theorem says that the space of holomorphic sections of
$\ca L_\lambda^n$ is isomorphic to $V_{n\lambda}$ as a $K$-module.
Hence, by definition, there exists $n>0$ such that $V_{n\lambda}$
contains an invariant vector if and only if $X\sst$ is nonempty.  By
Proposition \ref{proposition;semi} this is equivalent to the existence
of a $g\in G$ such that
\begin{equation}\label{equation;lambdainequal}
0\in f^*(w\inv\lambda+v\ca C)
\end{equation}
for all $(w,v)\in W\times W\rel$ satisfying
\begin{equation}\label{equation;transintersect}
gX^\circ_{wv}\cap\phi_v(\ti X)\ne\emptyset.
\end{equation}
Here we have used that $f^*(v\ca C)$ is dual to the cone spanned by
$\ti\Lambda_+\cap v\Lambda_+$, which follows from Lemma
\ref{lemma;dualcone}.  Since the semistable set is Zariski-open, by
Kleiman's transversality theorem (see e.g.\ Theorem 10.8 in
\cite[Ch. III]{hartshorne;algebraic-geometry;;1977}) we can perturb
$g$ to move $gX^\circ_{wv}$ into general position with respect to
$\phi_v(\ti X)$.  Then \eqref{equation;transintersect} amounts to
$\phi_v^!(c_{wv})\ne0$.  Now observe that the holomorphic embedding
$\phi_v\colon\ti X\to X$ is equal to the composition $v\inv\circ\phi$
(even though $v\inv$ is not holomorphic).  Hence
\begin{equation}\label{equation;dualnonzero}
0\ne\phi_v^!(c_{wv})=(v\inv\circ\phi)^!(c_{wv})=\phi^!(v\inv)^!(c_{wv})
=(-1)^{l(v)}\phi^!(vc_{wv}).
\end{equation}
Substituting $\lambda\to-w_0\lambda$ in
\eqref{equation;lambdainequal}, $w\to w_0w$ in
\eqref{equation;lambdainequal} and \eqref{equation;dualnonzero}, and
using that $V_{n\lambda}$ contains a $\ti K$-invariant vector if and
only if $V_{n\lambda}^*\cong V_{-nw_0\lambda}$ does, we obtain Theorem
\ref{theorem;invariant} for strictly dominant integral $\lambda$.

The case where $\lambda\in\Lambda^*_+$ is dominant but not strictly
dominant is very similar.  We briefly point out the main differences.
Instead of $X$ we use the partial flag variety $X_\s=K/K_\s\cong G/P$
with the $G$-equivariant line bundle $\ca
L_\lambda=G\times^P\C_\lambda$.  Here $\s$ is the face of $\t_+$ such
that $\lambda\in\s^\circ$, $P=P_\s$ is the parabolic subgroup of $G$
attached to $\s$, and $\C_\lambda$ is the obvious $P$-module defined
by $\lambda$.  For any $w\in W$ the closure of the $P$-orbit through
$\pi(w)$ is a single Schubert variety, namely $X_{w'}$, where $w'$ is
the longest element in the coset $W_\s w\subset W$.  The same
reasoning as above (using Proposition \ref{proposition;partialsemi}
instead of \ref{proposition;semi}) now gives that
$\lambda\in\Delta\bigl(T^*(\ti K\backslash K)\bigr)$ if and only if
\eqref{equation;lambdainequal} holds for all $w$ and $v$ satisfying
\begin{equation}\label{equation;pintersect}
\phi^!(vc_{w'v})\ne0.
\end{equation}
But if \eqref{equation;lambdainequal} holds for $w$, it holds for all
elements in the coset $W_\s w$, because $W_\s$ fixes $\lambda$.  This
finishes the proof of Theorem \ref{theorem;invariant}.

The proof of Theorem \ref{theorem;scalar} is almost exactly the same,
but uses Proposition \ref{proposition;maximal} instead of
\ref{proposition;semi}.  The details are left to the reader.

\section{Examples}\label{section;examples}

\subsection{Kostant's theorem}\label{subsection;kostant}

As a simple illustration of Theorem \ref{theorem;polytope} we consider
$\ti K=T$, the maximal torus of $K$.  The answer here is of course
well-known and can be obtained in many other ways.  Observe that
$\tilde T=T$, $\tilde W=\{1\}$, and every element of $\t$ is
``dominant'' with respect to $T$.  Therefore $W\rel=W$,
$\ti\t_v=v\t_+$, and the decomposition \eqref{equation;subdivision} is
simply the partition of $\t$ into Weyl chambers, $\t=\bigcup_{v\in
W}v\t_+$.  The dual cone to $v\t_+$ is $v\ca C$, where $\ca C$ is the
root cone of $K$.  The flag variety $\ti X$ is a point, so $\phi^*$ is
the trivial homomorphism $H^\bu(X,\Z)\to\Z$.  Since $\phi^*$ preserves
degree, we see that $\phi^*(v\sigma_{wv})\ne0$ if and only if
$v=w\inv$.  In other words, for every $\lambda\in\t^*_+$ the points
$\tlambda$ in the polytope $\Delta(\ca O_\lambda)$ are described by
the inequalities $w\tlambda\in\lambda-\ca C$, where $w$ ranges over
$W$.  This is equivalent to Kostant's result that $\Delta(\ca
O_\lambda)$ is the convex hull of the Weyl group orbit $W\lambda$.

\subsection{Klyachko's theorem}\label{subsection;klyachko}

In this section we apply our main theorem to the diagonal embedding of
a group into a number of copies of itself.  We deviate from our
standard notation and denote the small group by $K$ and the large
group $K\times K\times\dots\times K$ ($m$ times) by $K^m$.  The
canonical projection $f^*\colon(\t^*)^m\to\t^*$ is the addition map.
There is a canonical embedding of the small Weyl group $W$ into the
large Weyl group $W^m$, namely the diagonal embedding.  This tallies
with Example \ref{example;regular}, because every regular point in the
diagonal subgroup $K\hookrightarrow K^m$ is regular relative to $K^m$.
For the same reason, $W\rel=\{1\}$ and the partition
\eqref{equation;subdivision} has only one piece, namely $\t_1=\t_+$.
Its dual cone is $\ca C$, the root cone of $K$.  Moreover, the large
flag variety is $X^m$ and the embedding $\phi\colon X\to X^m$ is the
diagonal map.  The induced homomorphism $\phi^*\colon
H^\bu(X,\Z)^{\otimes m}\cong H^\bu(X^m,\Z)\to H^\bu(X,\Z)$ is the cup
product, and the Gysin map
$$
\phi^!\colon H_i(X^m,\Z)\to H_{i-(m-1)l(w_0)}(X,\Z)
$$
is the intersection product.  Similarly, for every face $\s$ of $\t_+$
the embedding \eqref{equation;embed} is simply the diagonal embedding
of the partial flag variety $X_\s$ into the $m$-fold product $X_\s^m$.
Now choose a collection of vectors $\chi_1$, $\chi_2,\dots$,
$\chi_n\in\t$ which span the cone $\t_+$.  Then Theorem
\ref{theorem;scalar} comes down to the following.

\begin{theorem}\label{theorem;klyachko}
Let $(\lambda_1,\lambda_2,\dots,\lambda_{m+1}) \in(\t^*_+)^{m+1}$.
Then $(\lambda_1,\lambda_2,\dots,\lambda_{m+1})
\in\Delta\bigl((T^*K)^m\bigr)$ if and only if
\begin{equation}\label{equation;add}
\sum_{l=1}^{m+1}\langle w_l\inv\lambda_l,\chi_k\rangle\ge0
\end{equation}
for $k=1$\upn, $2,\dots$\upn, $n$ and for all
$(w_1,w_2,\dots,w_{m+1})\in W^k\times W^k\times\dots\times W^k$ such
that $\sigma_{u_k}$ is contained in
$\sigma_{w_1}\cup\sigma_{w_2}\cup\dots\cup\sigma_{w_m}$.  Here $u_k\in
W^k$ denotes the shortest representative of $w_0w_{m+1}W_k$.
\end{theorem}

For $K=\U(n)$ this theorem was proved by Klyachko
\cite{klyachko;stable-bundles-hermitian}.  Also for $K=\U(n)$, Helmke
and Rosenthal \cite{helmke-rosenthal;eigenvalue} proved that the
inequalities \eqref{equation;add} are necessary, but not that they are
sufficient.  Inequalities for eigenvalues of sums of Hermitian
matrices have a long history going back to Weyl; see the cited papers
for examples.  The essential case of the theorem is of course when $K$
is semisimple; the identity component of the centre just contributes a
number of equalities, e.g.\ for the trace when $K=\U(n)$.  The cup
product condition is closely related to the Littlewood-Richardson
rules, which are discussed in Fulton
\cite{fulton;young,fulton;eigenvalues-sums}, Pragacz
\cite{pragacz;algebro-geometric} and Littelmann
\cite{littelmann;generalization}.

\begin{example}[$\G_2\times\G_2$]\label{example;g2}
Let us work out the case $K=\G_2$ and $m=2$.  Denote by
$S=\{\alpha_1,\alpha_2\}$ the simple roots of $\G_2$, where $\alpha_1$
is short and $\alpha_2$ is long.  Let $s_1=s_{\alpha_1}$ and
$s_2=s_{\alpha_2}$ be the corresponding simple reflections and
$\pi_1=2\alpha_1+\alpha_2$ and $\pi_2=3\alpha_1+2\alpha_2$ the
fundamental weights.  For $k=1$ or $2$, let $P_k$ be the parabolic
attached to $\pi_k$, $W_k$ the Weyl group, $W^k$ the set of shortest
representatives for $W/W_k$, and $X_k$ the associated ten-dimensional
``Grassmannian''.  Then
\begin{align*}
W_1 &=\{1,s_2\}, &W^1
&=\bigl\{1,s_1,s_2s_1,s_1s_2s_1,(s_2s_1)^2,s_1(s_2s_1)^2\bigr\},\\
W_2 &=\{1,s_1\}, &W^2
&=\bigl\{1,s_2,s_1s_2,s_2s_1s_2,(s_1s_2)^2,s_2(s_1s_2)^2\bigr\}.
\end{align*}
Using Chevalley's formula, Theorem
\ref{theorem;chevalley}\eqref{item;chevalley}, we can compute
multiplication tables for the cohomology of $X_k$ in terms of the
Schubert bases.  We find the relations
\begin{align*}
\sigma_{s_1}^2 &=\sigma_{s_2s_1},&\sigma_{s_1}^3
&=2\sigma_{s_1s_2s_1},&\sigma_{s_1}^4
&=2\sigma_{(s_2s_1)^2},&\sigma_{s_1}^5 &=2\sigma_{s_1(s_2s_1)^2}
&\text{for $k=1$},\\
\sigma_{s_2}^2 &=3\sigma_{s_1s_2},&\sigma_{s_2}^3
&=6\sigma_{s_2s_1s_2},&\sigma_{s_2}^4
&=18\sigma_{(s_1s_2)^2},&\sigma_{s_2}^5 &=18\sigma_{s_2(s_1s_2)^2}
&\text{for $k=2$},
\end{align*}
from which it is easy to derive all triples $(w_1,w_2,w_3)\in
W^k\times W^k\times W^k$ such that $\sigma_{u_k}$ is contained in
$\sigma_{w_1}\cup\sigma_{w_2}$.  For $k=1$ they are:
\begin{gather*}
\bigl(1,1,s_1(s_2s_1)^2\bigr),\quad\bigl(1,s_1,(s_2s_1)^2\bigr),\quad
(1,s_2s_1,s_1s_2s_1),\\
(s_1,s_1,s_1s_2s_1),\quad(s_1,s_2s_1,s_2s_1)
\end{gather*}
plus permutations of these triples.  For $k=2$ we obtain the same
list, but with $s_1$ and $s_2$ interchanged.  If we identify $\g_2^*$
with $\g_2$ by means of the invariant inner product for which
$\lVert\alpha_1\rVert^2=2$, then we can choose the $\chi_k$ to be the
fundamental weights $\pi_k$.  The inequalities for the polygons
$\Delta(\ca O_{\lambda_1}\times\ca O_{\lambda_2})$ are best written
out in coordinates relative to the basis $\{\pi_1,\pi_2\}$, because
then the positive chamber is given by the positive quadrant.  Writing
$\lambda_i=x_i\pi_1+y_i\pi_2$, we get for $k=1$
\begin{align*}
2x_1+3y_1+2x_2+3y_2-2x_3-3y_3 &\ge0\\
2x_1+3y_1+x_2+3y_2-x_3-3y_3 &\ge0\\
2x_1+3y_1+x_2-x_3 &\ge0\\
x_1+3y_1+x_2+3y_2-x_3 &\ge0\\
x_1+3y_1+x_2+x_3 &\ge0\\
\intertext{and for $k=2$}
x_1+2y_1+x_2+2y_2-x_3-2y_3 &\ge0\\
x_1+2y_1+x_2+y_2-x_3-y_3 &\ge0\\
x_1+2y_1+y_2-y_3 &\ge0\\
x_1+y_1+x_2+y_2-y_3 &\ge0\\
x_1+y_1+y_2+y_3 &\ge0,
\end{align*}
up to permutations of the variables.  The triples
$(s_1,s_2s_1,s_2s_1)$ and $(s_2,s_1s_2,s_1s_2)$ (the associated cycles
of which have intersection multiplicity $>1$) turn out to be
redundant.  The inequalities associated with the triples containing a
$1$ express the fact that
$$
\Delta(\ca O_{\lambda_1}\times\ca O_{\lambda_2})
\subset(\lambda_1+\hull W\lambda_2)\cap(\lambda_2+\hull W\lambda_1).
$$
Figure \ref{figure;g2} shows an example where this inclusion is
strict.  The dark shading indicates the moment polygon and the light
shading its Weyl group translates.  The dotted lines are the polygons
$\lambda_1+\hull W\lambda_2$ and $\lambda_2+\hull W\lambda_1$.

\begin{figure}
\setlength{\unitlength}{0.06mm}
$$
\begin{picture}(1300,1300)(-650,-650)
\thinlines
\dashline{10}(100,173.20508)(-100,-173.20508)
\dashline{10}(100,-173.20508)(-100,173.20508)
\dashline{10}(0,-173.20508)(0,173.20508)
\dashline{10}(-150,86.60254)(150,-86.60254)
\dashline{10}(-150,-86.60254)(150,86.60254)
\texture{cccccccc 0 0 0 cccccccc 0 0 0 cccccccc 0 0 0 cccccccc 0 0 0
	cccccccc 0 0 0 cccccccc 0 0 0 cccccccc 0 0 0 cccccccc 0 0 0 }
\shade\path(-50,86.60254)(-275,476.31397)(-225,562.91651)
(-150,606.21778)(-100,606.21778)(0,548.48276)(0,115.47005)(-50,86.60254)
\thicklines
\path(-50,86.60254)(-275,476.31397)(-225,562.91651)(-150,606.21778)
(-100,606.21778)(0,548.48276)(0,115.47005)(-50,86.60254)
\texture{c0c0c0c0 0 0 0 0 0 0 0 c0c0c0c0 0 0 0 0 0 0 0 c0c0c0c0 0 0 0
	0 0 0 0 c0c0c0c0 0 0 0 0 0 0 0 }
\thinlines
\shade\path(50,86.60254)(275,476.31397)(375,476.31397)(450,433.0127)
(475,389.71143)(475,274.24138)(100,57.735027)(50,86.60254)
\shade\path(100,57.735027)(475,274.24138)(575,216.50635)(600,173.20508)
(600,86.60254)(550,0)(100,0)(100,57.735027)
\shade\path(50,86.60254)(275,476.31397)(225,562.91651)(150,606.21778)
(100,606.21778)(0,548.48276)(0,115.47005)(50,86.60254)
\shade\path(-50,86.60254)(-275,476.31397)(-375,476.31397)(-450,433.0127)
(-475,389.71143)(-475,274.24138)(-100,57.735027)(-50,86.60254)
\shade\path(-100,57.735027)(-475,274.24138)(-575,216.50635)
(-600,173.20508)(-600,86.60254)(-550,0)(-100,0)(-100,57.735027)
\shade\path(50,-86.60254)(275,-476.31397)(225,-562.91651)
(150,-606.21778)(100,-606.21778)(0,-548.48276)(0,-115.47005)
(50,-86.60254)
\shade\path(50,-86.60254)(275,-476.31397)(375,-476.31397)(450,-433.0127)
(475,-389.71143)(475,-274.24138)(100,-57.735027)(50,-86.60254)
\shade\path(100,-57.735027)(475,-274.24138)(575,-216.50635)
(600,-173.20508)(600,-86.60254)(550,0)(100,0)(100,-57.735027)
\shade\path(-50,-86.60254)(-275,-476.31397)(-225,-562.91651)
(-150,-606.21778)(-100,-606.21778)(0,-548.48276)(0,-115.47005)
(-50,-86.60254)
\shade\path(-50,-86.60254)(-275,-476.31397)(-375,-476.31397)
(-450,-433.0127)(-475,-389.71143)(-475,-274.24138)(-100,-57.735027)
(-50,-86.60254)
\shade\path(-100,-57.735027)(-475,-274.24138)(-575,-216.50635)
(-600,-173.20508)(-600,-86.60254)(-550,0)(-100,0)(-100,-57.735027)
\put(-100,0){\circle*{15}}
\put(150,86.60254){\circle*{15}}
\put(0,0){\circle*{15}}
\put(-50,86.60254){\circle*{15}}
\put(0,173.20508){\circle*{15}}
\put(-150,606.21778){\circle*{15}}
\put(-125,303.10889){\circle*{15}}
\put(-25,303.10889){\circle*{15}}
\put(-160,20){\makebox(0,0)[lb]{\smash{$\alpha_1$}}}
\put(140,125){\makebox(0,0)[lb]{\smash{$\alpha_2$}}}
\put(-125,85){\makebox(0,0)[lb]{\smash{$\pi_1$}}}
\put(10,180){\makebox(0,0)[lb]{\smash{$\pi_2$}}}
\put(-160,625){\makebox(0,0)[lb]{\smash{$\lambda_1+\lambda_2$}}}
\put(-150,325){\makebox(0,0)[lb]{\smash{$\lambda_1$}}}
\put(-65,325){\makebox(0,0)[lb]{\smash{$\lambda_2$}}}
\dottedline{15}(-100,606.21778)(-150,606.21778)(-375,476.31397)
(-400,433.0127)(-400,173.20508)(-375,129.90381)(-150,0)(-100,0)
(125,129.90381)(150,173.20508)(150,433.0127)(125,476.31397)
(-100,606.21778)
\dottedline{15}(-150,606.21778)(-225,562.91651)(-350,346.41016)
(-350,259.80762)(-225,43.30127)(-150,0)(100,0)(175,43.30127)
(300,259.80762)(300,346.41016)(175,562.91651)(100,606.21778)
(-150,606.21778)
\end{picture}
$$
\caption{Moment ``rosette'' of $\ca O_{\lambda_1}\times\ca
O_{\lambda_2}$, where $K=\G_2$, $\lambda_1=\frac1{2}(5\pi_1+\pi_2)$,
and $\lambda_2=\frac1{2}(\pi_1+3\pi_2)$}
\label{figure;g2}
\end{figure}
\end{example}

\subsection{$\lie{sl}(2)$-triples}\label{subsection;sl2}

Our next example is a general homomorphism with finite kernel
$f\colon\ti K\to K$, where $\ti K=\SU(2)$ and $K$ is semisimple.  Up
to conjugacy such homomorphisms are in one-to-one correspondence with
embeddings of $\ti\g=\lie{sl}(2,\C)$ into $\g$, that is to say,
triples $(\ti h,\ti e,\ti f)$ of vectors in $\g$ which satisfy
$$
[\ti h,\ti e]=2\ti e,\qquad[\ti h,\ti f] =-2\ti f,\qquad[\ti e,\ti
f]=\ti h.
$$
Indeed, because $\SL(2,\C)$ is simply connected any such embedding
lifts to a homomorphism $\SL(2,\C)\to G$, which can be conjugated to a
homomorphism that maps $\SU(2)$ into $K$, since $K$ is maximal compact
in $G$.  In other words, we may assume that $\ti h\in\t$, $\ti
e\in\lie n$ and $\ti f\in\theta\lie n$.  After a further conjugation
with an element of $W$ we may even assume that $\ti h\in\t_+$.  Then
the chambers $\ti\t_+$ and $\t_+$ are compatible and $W\rel=\{1\}$.

Let $\{\,\alpha^*\mid\alpha\in S\,\}$ be the set of fundamental
coweights, i.e.\ the basis of $\t$ which is dual to $S\subset\t^*$.
Then we can write
\begin{equation}\label{equation;h}
\ti h=\sum_{\alpha\in S\setminus\bar S}d_\alpha\alpha^*,
\end{equation}
where all $d_\alpha>0$ and $\bar S\subset S$ is the set of simple
roots orthogonal to $\ti h$.  According to Dynkin's classification of
$\lie{sl}(2)$-triples (see \cite{dynkin;semisimple-subalgebras}),
$d_\alpha=1$ or $2$ for all $\alpha\not\in\bar S$.  The set $\bar S$
is a base of the root system $\bar R$ of $\bar K=\ca Z_K(\ti T)$, and
$\barW$ is generated by the simple reflections $s_\alpha$ with
$\alpha\in\bar S$.  Let $j\colon\tW=\{1,\tilde w_0\}\to W$ be the
inclusion map determined by the choice of dominant chambers $\ti\t_+$
and $\t_+$.  This map is easy to describe explicitly.

\begin{lemma}\label{lemma;image}
\begin{enumerate}
\item\label{item;w0}
$w_0\ti h=-\ti h$.
\item\label{item;short}
$j(\ti w_0)=w_0\bar w_0$.
\item\label{item;involution}
$w_0\bar w_0$ is an involution.
\end{enumerate}
\end{lemma}

\begin{proof}
The $\tW$-equivariance of the embedding $f_*\colon\ti\t\to\t$
implies that $j(\ti w_0)\ti h=-\ti h$, which is in the antidominant
chamber $-\t_+$.  On the other hand $w_0\ti h$ is also in $-\t_+$ and
hence $w_0\ti h=-\ti h$.  This proves \eqref{item;w0}.

To prove \eqref{item;short} it suffices to show that $w_0\bar w_0$ is
in $\ca N_W(\bar S)$ and maps $\ti h$ to $-\ti h$.  It follows from
\eqref{item;w0} that $w_0$ preserves $\bar R$, and therefore sends
$\bar R_+$ to $\bar R_-$.  Consequently $w_0\bar w_0$ preserves $\bar
S$.  Since $\barW$ fixes $\ti h$, $w_0\bar w_0\ti h=w_0\ti h=-\ti h$
by \eqref{item;w0}.

\eqref{item;involution} follows immediately from \eqref{item;short}.
\end{proof}

In order to write the inequalities for the ``polytopes'' $\Delta(\ca
O_\lambda)$ we identify $\ti\t^*$ with $\R$ by sending the positive
root $\ti\alpha$ to $2$.  Dually, this corresponds to sending the
basis element $\ti h\in\ti\t$ to $1$.  The projection
$f^*(\lambda)\in\ti\t^*$ of any functional $\lambda\in\t^*$ then gets
identified with the number $\lambda(\ti h)$.

\begin{proposition}
Let $\lambda\in\t^*_+$.  Then $\tlambda\in\Delta(\ca O_\lambda)$ if
and only if $\tlambda\ge0$ and
$$
-\lambda(\ti h)+\max_{\alpha\in S\setminus\bar S}
d_\alpha\lambda(\alpha\spcheck)\le\tlambda\le\lambda(\ti h).
$$
\end{proposition}

\begin{proof}
Following Theorem \ref{theorem;polytope} we determine all $(\ti
w,w)\in\tW\times W$ such that $\phi^*(\sigma_w)(\ti c_{\ti w})\ne0$,
i.e.\ $\ti\sigma_{\ti w}$ is contained in $\phi^*(\sigma_w)$.  Since
$\ti X=\C P^1$, $H^{2l(\ti w)}(\ti X,\Z)$ vanishes for $l(\ti w)>1$.
Moreover, $\phi^*$ preserves degree, so that we need only consider
Weyl group elements of length $\le1$.

For $\ti w=w=1$ we find $\ti\sigma_{\ti w}=\sigma_w=1$ and get the
inequality $\tlambda\le\lambda(\ti h)$.

If $l(\ti w)=l(w)=1$ then $\ti w=\ti w_0$ and $w=s_\alpha$ for some
$\alpha\in S$.  Using Theorem \ref{theorem;chevalley} we find that
$\phi^*(\sigma_{s_\alpha})=f^*(\pi_\alpha)\ti\sigma_{\ti
w_0}=\pi_\alpha(\ti h)\ti\sigma_{\ti w_0}$, where $\pi_\alpha$ denotes
the fundamental weight corresponding to $\alpha$.  Therefore
$\ti\sigma_{\ti w_0}$ is contained in $\phi^*(\sigma_{s_\alpha})$ if
and only if $\pi_\alpha(\ti h)\ne0$.  From \eqref{equation;h} we
obtain
$$
\pi_\alpha(\ti h)=\sum_{\beta\in S\setminus\bar
S}d_\beta\pi_\alpha(\beta^*).
$$
Since the angle between any two fundamental weights $\pi_\alpha$ and
$\pi_\beta$ is acute, we have $\pi_\alpha(\beta^*)\ge0$ for all
$\alpha$ and $\beta$, and therefore $\pi_\alpha(\ti h)\ne0$ if and
only if $\pi_\alpha(\beta^*)\ne0$ for some $\beta\in S\setminus\bar
S$.  This is certainly the case if $\alpha\in S\setminus\bar S$,
because then we can take $\beta=\alpha$, but may also happen if
$\alpha\in\bar S$.  In any case, the inequality corresponding to $(\ti
w_0,s_\alpha)$ is $\ti w_0\tlambda\le f^*(s_\alpha\ti h)
=\lambda(s_\alpha\lambda)$, or
$$
\tlambda\ge-\lambda(s_\alpha\ti h) =-\lambda(\ti
h)+\lambda(\alpha\spcheck)\alpha(\ti h).
$$
Here we have used that $\ti w_0\tlambda=-\tlambda$ and
$s_\alpha\lambda =\lambda-\lambda(\alpha\spcheck)\alpha$.  Moreover,
$$
\alpha(\ti h)=\sum_{\beta\in S\setminus\bar S}d_\beta\alpha(\beta^*) =
\begin{cases}
d_\alpha &\text{if $\alpha\not\in\bar S$}\\
0 &\text{if $\alpha\in\bar S$.}
\end{cases}
$$
Thus, for $\alpha\in\bar S$ we obtain the inequality
$\tlambda\ge-\lambda(\ti h)$, which is vacuous, and for
$\alpha\not\in\bar S$ we obtain $\tlambda\ge-\lambda(\ti
h)+\lambda(\alpha\spcheck)d_\alpha$.
\end{proof}

\begin{example}\label{example;principal}
There exists a triple $(\ti h,\ti e,\ti f)$, known as the
\emph{principal} triple, for which $\bar S=\emptyset$ and $d_\alpha=2$
for all $\alpha$.  It was shown by Dynkin that principal triples are
unique up to conjugation.  In this case we get the inequalities
$-\lambda(\ti h)+2\max_{\alpha\in S}\lambda(\alpha\spcheck)
\le\tlambda\le\lambda(\ti h)$.  For instance, if $K$ is the product of
$m$ copies of $\SU(2)$, then the principal $\SU(2)$ is the diagonal
subgroup.  Here $\lambda$ can be represented as an $m$-tuple
$(\tlambda_1,\tlambda_2,\dots,\tlambda_m)$ and the inequalities are
$$
-\tlambda_1-\tlambda_2-\dots-\tlambda_m+2\max_i\tlambda_i
\le\tlambda\le\tlambda_1+\tlambda_2+\dots+\tlambda_m.
$$
This can of course also be regarded as a special case of Theorem
\ref{theorem;klyachko}.
\end{example}

\subsection{A maximal rank
subgroup}\label{subsection;maximal}

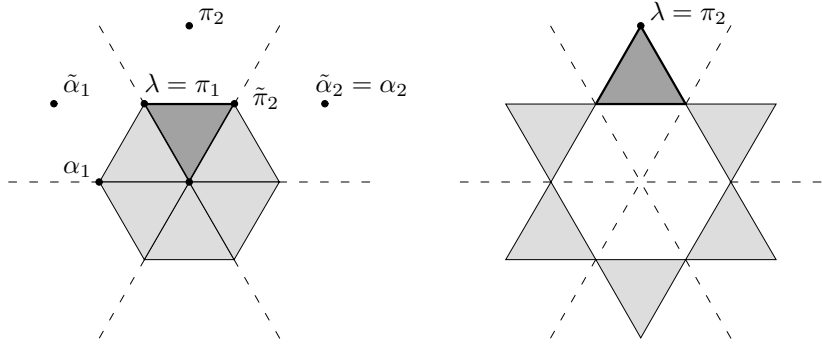
\begin{figure}
\setlength{\unitlength}{0.12mm}
$$
\begin{picture}(500,500)(-250,-250)
\thinlines
\dashline{10}(100,173.20508)(-100,-173.20508)
\dashline{10}(-200,0)(200,0)
\dashline{10}(100,-173.20508)(-100,173.20508)
\texture{cccccccc 0 0 0 cccccccc 0 0 0 cccccccc 0 0 0 cccccccc 0 0 0
	cccccccc 0 0 0 cccccccc 0 0 0 cccccccc 0 0 0 cccccccc 0 0 0 }
\shade\path(-50,86.60254)(50,86.60254)(0,0)(-50,86.60254)
\thicklines
\path(-50,86.60254)(50,86.60254)(0,0)(-50,86.60254)
\texture{c0c0c0c0 0 0 0 0 0 0 0 c0c0c0c0 0 0 0 0 0 0 0 c0c0c0c0 0 0 0
	0 0 0 0 c0c0c0c0 0 0 0 0 0 0 0 }
\thinlines
\shade\path(-50,-86.60254)(50,-86.60254)(0,0)(-50,-86.60254)
\shade\path(-50,86.60254)(0,0)(-100,0)(-50,86.60254)
\shade\path(50,86.60254)(0,0)(100,0)(50,86.60254)
\shade\path(-50,-86.60254)(0,0)(-100,0)(-50,-86.60254)
\shade\path(50,-86.60254)(0,0)(100,0)(50,-86.60254)
\put(-100,0){\circle*{7}}
\put(-150,86.60254){\circle*{7}}
\put(150,86.60254){\circle*{7}}
\put(0,0){\circle*{7}}
\put(-50,86.60254){\circle*{7}}
\put(50,86.60254){\circle*{7}}
\put(0,173.20508){\circle*{7}}
\put(-140,10){\makebox(0,0)[lb]{\smash{$\alpha_1$}}}
\put(-140,100){\makebox(0,0)[lb]{\smash{$\ti\alpha_1$}}}
\put(140,100){\makebox(0,0)[lb]{\smash{$\ti\alpha_2=\alpha_2$}}}
\put(-50,100){\makebox(0,0)[lb]{\smash{$\lambda=\pi_1$}}}
\put(70,85){\makebox(0,0)[lb]{\smash{$\ti\pi_2$}}}
\put(10,180){\makebox(0,0)[lb]{\smash{$\pi_2$}}}
\end{picture}
\begin{picture}(500,500)(-250,-250)
\thinlines
\dashline{10}(100,173.20508)(-100,-173.20508)
\dashline{10}(-200,0)(200,0)
\dashline{10}(100,-173.20508)(-100,173.20508)
\texture{cccccccc 0 0 0 cccccccc 0 0 0 cccccccc 0 0 0 cccccccc 0 0 0
	cccccccc 0 0 0 cccccccc 0 0 0 cccccccc 0 0 0 cccccccc 0 0 0 }
\shade\path(-50,86.60254)(50,86.60254)(0,173.20508)(-50,86.60254)
\thicklines
\path(-50,86.60254)(50,86.60254)(0,173.20508)(-50,86.60254)
\texture{c0c0c0c0 0 0 0 0 0 0 0 c0c0c0c0 0 0 0 0 0 0 0 c0c0c0c0 0 0 0
	0 0 0 0 c0c0c0c0 0 0 0 0 0 0 0 }
\thinlines
\shade\path(-50,-86.60254)(50,-86.60254)(0,-173.20508)(-50,-86.60254)
\shade\path(-50,86.60254)(-150,86.60254)(-100,0)(-50,86.60254)
\shade\path(50,86.60254)(150,86.60254)(100,0)(50,86.60254)
\shade\path(-50,-86.60254)(-150,-86.60254)(-100,0)(-50,-86.60254)
\shade\path(50,-86.60254)(150,-86.60254)(100,0)(50,-86.60254)
%
\put(0,173.20508){\circle*{7}}
\put(10,180){\makebox(0,0)[lb]{\smash{$\lambda=\pi_2$}}}
\end{picture}
$$
\caption{Moment rosettes of $\G_2$-orbits $\ca O_\lambda$ w.r.t.\
$\SU(3)$-action, where $\lambda=\pi_1$, resp.\ $\lambda=\pi_2$
(complexified adjoint representation)}
\label{figure;g2a2}
\end{figure}

Let $K=\G_2$.  The long roots form a subsystem of the root system of
$K$ and the associated subgroup $\ti K$ is isomorphic to $\SU(3)$.
Using the notation of Example \ref{example;g2}, we denote the simple
roots of $K$ by $\alpha_1$ and $\alpha_2$ and its fundamental weights
by $\pi_1$ and $\pi_2$.  Writing $\ti\alpha_1$ and $\ti\alpha_2$ for
the simple roots of $\ti K$, and $\ti\pi_1$ and $\ti\pi_2$ for its
fundamental weights, we have
$\alpha_1=\frac1{3}(\ti\alpha_1-\ti\alpha_2)$, $\alpha_2=\ti\alpha_2$,
$\pi_1=\ti\pi_1$ and $\pi_2=\ti\pi_1+\ti\pi_2$.  The matrix of $f^*$
relative to the fundamental weights is therefore
\begin{equation}\label{equation;fstar}
f^*=\begin{pmatrix}1&1\\0&1\end{pmatrix}.
\end{equation}
The dominant chamber of $K$ consists of two chambers of $\ti K$, and
the relative Weyl set is $\{1,s_1\}$.  According to Theorem
\ref{theorem;scalar} we need to consider three one-parameter subgroups
$\ti\chi_1$, $\ti\chi_2$, $\ti\chi_3$, and corresponding embeddings
$\phi_1\colon\ti G/\ti P_1\to G/P_1$, $\phi_2\colon\ti G/\ti B\to
G/P_2$ and $\phi_3\colon\ti G/\ti P_2\to G/P_1$.  However, the
inequalities coming from $\phi_3$ are dual to those of $\phi_1$, so we
need to consider only $\phi_1$ and $\phi_2$.  The cohomology of the
Grassmannians of $K$ was computed in Example \ref{example;g2}.  The
calculation for $\ti K$ is very similar.  By means of
\eqref{equation;fstar} and Theorem \ref{theorem;chevalley} we can then
calculate all pairs $(\ti w,w)\in\tW^k\times W^k$ such that
$\ti\sigma_{\ti w}$ is contained in $\phi_k^*(\sigma_w)$.  For $k=1$
we find
$$
(1,1),\qquad(\ti s_1,s_1),\qquad(\ti s_2\ti s_1,s_2s_1),
$$
and for $k=2$,
$$
(1,1),\quad(\ti s_1,s_2),\quad(\ti s_2,s_2),\quad(\ti s_1\ti
s_2,s_1s_2),\quad(\ti s_2\ti s_1,s_1s_2),\quad(\ti w_0,s_2s_1s_2).
$$
In fact, for all these pairs $\ti\sigma_{\ti w}$ occurs with
multiplicity $1$ in $\phi_k^*(\sigma_w)$.  Introducing coordinates
$\tlambda=\ti x\ti\pi_1+\ti y\ti\pi_2$ and $\lambda=x\pi_1+y\pi_2$, we
can then write the inequalities $\langle\ti
w\inv\tlambda,\ti\chi_k\rangle\le\langle w\inv\lambda,\chi_k\rangle$
as follows.  It turns out that the inequalities for $k=1$ are all
redundant, as are their dual inequalities (which are obtained by
interchanging $\ti x$ and $\ti y$).  For $k=2$ the pairs $(\ti s_1\ti
s_2,s_1s_2)$ and $(\ti s_2\ti s_1,s_1s_2)$ lead to redundant
inequalities.  The remaining ones are
$$
\ti x+\ti y\le x+2y,\quad\ti x\le x+y,\quad\ti y\le x+y,\quad\ti x+\ti
y\ge y.
$$
This set of inequalities is self-dual (stable under interchanging $\ti
x$ and $\ti y$).  In addition to the obvious inequalities $\ti x\ge0$,
$\ti y\ge0$, $x\ge0$, $y\ge0$, these completely describe the set of
all $(\tlambda,\lambda)\in\ti\t^*_+\times\t^*_+$ such that
$\tlambda\in\Delta(\ca O_\lambda)$.  See Figure \ref{figure;g2a2} for
examples.

\appendix

\section{Flag varieties}\label{section;flag}

\subsection{Schubert cells}\label{subsection;cell}

Let $G=K^\C$ be the complexification of $K$ and $H=T^\C$ the
complexified maximal torus.  Let $R\subset\t^*$ be the root system of
$K$ and $R_+$ the set of positive roots.  Let $\lie n$ be the
nilpotent subalgebra of $\g$ spanned by the positive root spaces and
$N=\exp\lie n$ the corresponding maximal unipotent subgroup of $G$.
Furthermore, let $\lie b$ be the Borel subalgebra $\h\oplus\lie n$ and
$B=\exp\lie b$ the corresponding Borel subgroup.  The canonical map
$$
\tau\colon X=K/T\to G/B
$$
defined by $\tau(kT/T)=kB/B$ is a $K$-equivariant diffeomorphism.  The
complex homogeneous space $G/B$ decomposes into Bruhat cells
$X^\circ_w=BwB/B$, where $w\in W$.  The \emph{Schubert variety} $X_w$
is the closure of $X^\circ_w$; the \emph{Schubert class} is the
fundamental class $[X_w]\in H_{2l(w)}(G/B,\Z)$, where $l(w)$ is the
length of $w$.  Define
$$
c_w =\tau_*\inv\bigl([X_w]\bigr);
$$
these classes form the \emph{Schubert basis} of the homology group
$H_\bu(X,\Z)$.

\begin{remark}\label{remark;basis}
The Schubert basis of $H_\bu(X,\Z)$ depends on the choice of the set
of positive roots.  Let us work out the formula for a change of basis.
Any set of positive roots can be written as $uR_+$ for a unique $u\in
W$, and the Borel associated with $uR_+$ is $B^u=uBu\inv$.  Let us
write $(X^u_w)^\circ$ for the Schubert cell $B^uwB^u/B^u$ in $G/B^u$,
$\tau_u$ for the canonical $K$-equivariant diffeomorphism $X\to
G/B^u$, and $c^u_w$ for the homology class
$(\tau_u)_*\inv\bigl([X^u_w]\bigr)$.  Consider the commutative diagram
$$
\begin{CD}
X@>\tau>>G/B\\
@V{u}VV@VV{\psi_u}V\\
X@>{\tau_u}>>G/B^u,
\end{CD}
$$
where the vertical arrow on the left denotes the action of $u$ on $X$
and the map $\psi_u$ is the $G$-equivariant holomorphic map
$\psi_u(gB/B)=gu\inv B^u/B^u$.  It is clear that $\psi_u(X_w^\circ)
=u\inv(X^u_{uwu\inv})^\circ$, so we see that
$$
(\tau_u)_*(uc_w) =(\psi_u)_*\tau_*(c_w) =(\psi_u)_*\bigl([X_w]\bigr)
=[X^u_{uwu\inv}] =(\tau_u)_*c^u_{uwu\inv}.
$$
We conclude that the Schubert basis relative to the set of positive
roots $uR_+$ is given by $c^u_w=uc_{u\inv wu}$.
\end{remark}

If $P$ is any parabolic subgroup of $G$, we have the partial flag
variety $G/P$ and a $K$-equivariant diffeomorphism $\tau\colon
K/(P\cap K)\to G/P$.  We call $P$ \emph{standard} if $P\supset B$;
then $P=G_\s N$, where $G_\s$ is the complexification of $K_\s=\ca
Z_K(\s)$, the centralizer of a face $\s\subset\t_+$.  In this case we
write $P=P_\s$ and $K/K_\s=X_\s$.  Every parabolic is of the form
\begin{equation}\label{equation;parabolic}
P_\chi=\bigl\{\,p\in G\bigm|
\lim_{t\to0}\chi(t)\,p\,\chi(t)\inv\quad\text{exists}\,\bigr\}
\end{equation}
for some algebraic one-parameter subgroup $\chi\in\Hom(\C^\times,G)$,
and we have $P_\chi=P_\s$ if and only if (the infinitesimal generator
of) $\chi$ is in $\s^\circ$.  The \emph{Schubert cells} in $G/P$ are
the sets $(X_\s)_w^\circ=BwP/P$; their closures are the \emph{Schubert
varieties} $(X_\s)_w$, and the classes $c^\s_w
=\tau_*\inv\bigl([(X_\s)_w]\bigr)$ form the \emph{Schubert basis} of
$H_\bu(X_\s,\Z)$.  Here $w$ ranges over the subset $W^\s$ of $W$.

\subsection{Cohomology}\label{subsection;cohomology}

Consider the map $H^\bu(X,\Z)\to H^\bu(\ti X,\Z)$ induced by the
embedding $\phi\colon\ti X\to X$.  Note that
$X=[K,K]\big/\bigl(T\cap[K,K]\bigr)$ and that $f$ maps $[\ti K,\ti K]$
into $[K,K]$.  Hence $\phi$ can also be viewed as the embedding
induced by the homomorphism $[\ti K,\ti K]\to[K,K]$.  For the purpose
of this discussion we may therefore assume both $K$ and $\ti K$ to be
semisimple.

Define a map $\Theta\colon\Lambda^*\to H^2(X,\Z)$ by
$\Theta(\lambda)=c_1(\ca L_\lambda)$, the first Chern class of the
homogeneous line bundle with weight $\lambda$.  See
\cite{bernstein-gelfand-gelfand;schubert} or
\cite{demazure;desingularisation} for the theorems quoted below.

\begin{theorem}[Chevalley]\label{theorem;chevalley}
\begin{enumerate}
\item\label{item;isomorphism}
$\Theta$ is an isomorphism.
\item\label{item;fundamental}
$\Theta(\pi_\alpha)=\sigma_{s_\alpha}$ for all simple roots $\alpha$.
\item\label{item;chevalley}
For all weights $\lambda$
$$
\Theta(\lambda)\cup\sigma_w =\sum_{\substack{\beta\in R_+\\l(ws_\beta)
=l(w)+1}} \lambda(\beta\spcheck)\sigma_{ws_\beta}.
$$
\end{enumerate}
\end{theorem}

Parts \eqref{item;isomorphism} and \eqref{item;fundamental} lead to
the following simple description of the map $\phi^*$ restricted to
$H^2$.  It is clear that the diagram
$$
\begin{CD}
\Lambda^*@>{f^*}>>\ti\Lambda^*\\
@V{\Theta}VV@VV{\ti\Theta}V\\
H^2(X,\Z)@>{\phi^*}>>H^2(\ti X,\Z)
\end{CD}
$$
commutes.  In degree $2$, the matrix of $\phi^*$ relative to the
Schubert bases is therefore the same as the matrix of $f^*$ relative
to the bases of fundamental weights.

Now consider the graded algebra $S(\t^*)$ of polynomial functions on
$\t$ (in which elements of $\t^*$ are defined to be of degree $2$) and
the ideal $J$ generated by the $W$-invariant polynomials of positive
degree.  Recall the following well-known result.

\begin{theorem}[Borel]\label{theorem;borel}
The map $\Theta$ extends to a surjective homomorphism of graded
algebras $S(\t^*)\to H^*(X,\R)$\upn, whose kernel is equal to $J$.
\end{theorem}

Observe that $f^*(J)\subset\ti J$, so that the map $\phi^*$ can be
alternatively described (at least over $\R$) as the homomorphism
induced by the restriction map $S(\t^*)\to S(\ti\t^*)$.  Thus the
degree-$2$ Schubert classes $\sigma_{s_\alpha}$ generate the
cohomology of $X$ over $\R$ (though not over $\Z$), and by using
Theorem \ref{theorem;chevalley}\eqref{item;chevalley} one can in
principle calculate the matrix of $\phi^*$ in higher degrees.  This
can be laborious in practice, however.

A convenient method for calculating the Weyl group action on the
cohomology is provided by the difference operators $D_\alpha\colon
S(\t^*)\to S(\t^*)$, which are defined for $\alpha\in R$ by
$$
D_\alpha(p)=\frac{p-s_\alpha(p)}{\alpha},
$$
where $p\in S(\t^*)$.  It is clear that $D_\alpha(J)=\{0\}$.  Hence,
by Theorem \ref{theorem;borel}, $D_\alpha$ induces a linear operator
of degree $-2$ on $H^*(X,\R)$, which will also be denoted by
$D_\alpha$.  Note that $D_\alpha^2=0$ for all $\alpha$.

\begin{theorem}[Bernstein et al.\
\cite{bernstein-gelfand-gelfand;schubert}, Demazure
\cite{demazure;desingularisation}]
\label{theorem;dbgg}
For all simple roots $\alpha$ and all $w\in W$ we have
$$
D_\alpha(\sigma_w)=
\begin{cases}
0 &\text{if $ws_\alpha\succ w$}\\
\sigma_{ws_\alpha} &\text{if $ws_\alpha\prec w$.}
\end{cases}
$$
\end{theorem}

Thus $s_\alpha\sigma_w=\sigma_w-\Theta(\alpha)D_\alpha(\sigma_w)$,
which is equal to $\sigma_w$ if $ws_\alpha\succ w$ and to
$\sigma_w-\Theta(\alpha)\cup\sigma_{ws_\alpha}$ if $ws_\alpha\prec w$.

\section{Notation}\label{section;notation}

\begin{tabbing}
\indent \= $M_\mu$; $M_0=M\qu G$ \= \kill 
\> $K$; $T$; $C$ \> compact connected Lie group; maximal torus;
centre\\
\> $G$; $H$ \> complexification of $K$; resp.\ $T$\\
\> $W$; $w_0$ \> Weyl group of $K$ w.r.t.\ $T$; longest Weyl group
element\\
\> $R$; $R_+$ \> root system; positive roots\\
\> $S$; $\ca C$ \>  simple roots; cone spanned by $S$\\
\> $N$; $B$ \> maximal unipotent subgroup of $G$; Borel subgroup
$HN$\\
\> $\ca N_I(J)$; $\ca Z_I(J)$ \> normalizer; resp.\ centralizer of
group $J$ in a group $I$\\
\> $X$; $\pi$ \> full flag variety $K/T\cong G/B$; projection $G\to
X$\\
\> $\t_+$; $\t^*_+$ \> positive Weyl chamber in $\lie t$; positive
Weyl chamber in $\lie t^*$\\
\> $\Lambda$; $\Lambda^*$ \> integral lattice in $\lie t$; weight
lattice $\Hom_\Z(\Lambda,\Z)$\\
\> $\Lambda_+$; $\Lambda^*_+$ \> dominant one-parameter subgroup;
dominant characters\\
\> $\ca O_\lambda$ \> coadjoint orbit through $\lambda\in\k^*$\\
\> $f$; $f^*$ \> homomorphism $\ti K\to K$; dual projection
$\k^*\to\ti\k^*$\\
\> $X^\circ_w$; $X_w$ \> Bruhat cell $BwB/B\subset X$; its closure\\
\> $c_w$; $\sigma_w$ \> Schubert cycle $[X_w]$ in $H_{2l(w)}(X,\Z)$;
cocycle in $H^{2l(w)}(X,\Z)$\\
\> $\ca P$ \> Poincar\a'e duality $H^i(X,\Z)\to H_{2l(w_0)-i}(X,\Z)$\\
\> $\phi$; $\phi^!$\> embedding $\ti X\to X$; wrong-way map $\ti\ca
P\circ\phi^*\circ\ca P\inv$\\
\> $\s$; $\s^\circ$ \> face of $\t_+$; its relative interior\\
\> $K_\s$; $R_\s$ \> centralizer of $\s$; its root system\\
\> $W_\s$; $W^\s$ \> Weyl group of $K_\s$; set of shortest
representatives of $W/W_\s$\\
\> $P_\s$; $X_\s$ \> parabolic $G_\s N$; partial flag variety
$K/K_\s\cong G/P_\s$\\
\> $\pi_\s$; $\phi_{\ti\s}$ \> projection $G\to X_\s$; inclusion $\ti
X_{\ti\s}\to X_\s$ (if $\ti K_{\ti\s}\subset K_\s$)\\
\> $\ca L_\lambda$ \> line bundle $G\times^{P_\s}\C_\lambda$ on $X_\s$
(for $\lambda\in\s^\circ\cap\Lambda^*$)\\
\> $P_\chi$ \> parabolic associated with one-parameter subgroup $\chi$
of $G$\\
\> $(X_\s)^\circ_w$; $(X_\s)_w$ \> Bruhat cell $BwP/P\subset X_\s$;
its closure\\
\> $c^\s_w$ \> Schubert cycle $[(X_\s)_w]$ in $H_{2l(w)}(X_\s,\Z)$
(for $w\in W^\s$)\\
\> $\ti\t_w$; $W\rel$ \> cone $w\t_+\cap\ti\t$ (for $w\in W$);
relative Weyl set
\end{tabbing}

Similar conventions are in force for the groups $\ti K$ and $\bar K$,
i.e.\ $\ti C$ and $\bar C$ denote the centres of $\ti K$ and $\bar K$,
$\ti G$ and $\bar G$ their complexifications, etc.  We assume
furthermore that $f(\ti T)\subset T$ and $f(\ti B)\subset B$.  Here
$\ti K$ is an arbitrary compact connected Lie group, $f\colon\ti K\to
K$ is a homomorphism with finite kernel, and $\bar K=\ca
Z_K\bigl(f(\ti T)\bigr)$.

\subsection*{Weights, roots}

We identify the character group $\Hom(T,S^1)$ with the weight lattice
$\Lambda^*$ by mapping a character $\lambda$ to the weight
$d\lambda(1)/2\pi i$.  We identify $\lie c^*$ with the annihilator of
$[\k,\k]$ in $\k^*$ and $[\k,\k]^*$ with the annihilator of $\lie c$,
so that we have a direct sum decomposition $\k^*=\lie
c^*\oplus[\k,\k]^*$.

The \emph{roots} are the weights of the complexified adjoint
representation $\g$.  Over $\R$ they span the subspace
$\t^*\cap[\k,\k]^*$ of $\t^*$.  Given a base $S$ of $R$, the set of
\emph{fundamental coweights} is the basis $\{\,\alpha^*\mid\alpha\in
S\,\}$ of $\t\cap[\k,\k]$ which is dual to $S$.  The \emph{dual root
system} $R\spcheck$ is the root system in $\t\cap[\k,\k]$ consisting
of the \emph{dual roots} or \emph{coroots} $\alpha\spcheck$, which are
determined by $s_\alpha\lambda
=\lambda-\lambda(\alpha\spcheck)\alpha$.  The set of \emph{fundamental
weights} is the basis $\{\,\pi_\alpha\mid\alpha\in S\,\}$ of
$\t^*\cap[\k,\k]^*$ which is dual to $S\spcheck$.

If we identify $\t$ with $\t^*$ via a $W$-invariant inner product
$\langle\cdot,\cdot\rangle$, then $\alpha\spcheck
=2\alpha/\langle\alpha,\alpha\rangle$ and $\alpha^*
=2\pi_\alpha/\langle\alpha,\alpha\rangle$.

The Weyl chamber $\t_+$ is equal to $\lie
c\times\cone\{\,\alpha^*\mid\alpha\in S\,\}$, whereas $\t^*_+$ is
equal to $\lie c^*\times\cone\{\,\pi_\alpha\mid\alpha\in S\,\}$.
(Here $\cone(X)$ denotes the set of all nonnegative linear
combinations of elements in $X$.)  We define the \emph{root cone} to
be the cone $\ca C$ spanned by $R_+$; it is dual to the cone $\t_+$ in
the sense that $\lambda\in\ca C$ if and only if $\lambda(\xi)\ge0$ for
all $\xi\in\t_+$.


\bibliographystyle{amsplain}

\providecommand{\bysame}{\leavevmode\hbox to3em{\hrulefill}\thinspace}


\end{document}